\documentclass{article}

\usepackage[T1]{fontenc}
\usepackage[utf8]{inputenc}
\usepackage{soul}
\usepackage{import}
\usepackage{hyperref}

\usepackage{graphicx}

\usepackage{amsmath,amsthm,amssymb,amsfonts,amscd, amsbsy,mathtools}
\usepackage{wasysym}
\usepackage{mathrsfs}

\usepackage{empheq}

\usepackage{ upgreek }

\usepackage{pdfpages}

\usepackage{bussproofs}

\usepackage{subcaption}

\everymath{\displaystyle}
\DeclareMathAlphabet{\mathcal}{OMS}{cmsy}{m}{n}

\usepackage{xcolor}

\usepackage{tikz}

\usepackage[margin=2cm]{geometry}

\usetikzlibrary{arrows.meta}

\definecolor{curveblue}{RGB}{39,91,169}
\definecolor{curvered}{RGB}{207,57,63}

\usepackage{faktor}

\usepackage{subcaption}
\usepackage{float}

\usepackage{amsmath}
\usepackage{graphicx}
\newtheorem{remark}{Remark}[section]
\newtheorem{prop}{Proposition}[section]
\newtheorem{theorem}{Theorem}[section]
\newtheorem{lemma}{Lemma}[section]
\newtheorem{cor}{Corollary}[section]

\title{Weighted minimal surfaces in $\mathbb{H}^2 \times \mathbb{R}$ invariant under hyperbolic translations}
\author{
    \fontsize{12}{14}\selectfont
    Katianny Rosa Freitas \textsuperscript{1*}, João Paulo dos Santos \textsuperscript{2}, \fontsize{12}{14}\selectfont
    António Martínez \textsuperscript{3}\\[1ex]
    \fontsize{12}{12}\selectfont
    \textsuperscript{1, 2}Departamento de Matemática, Universidade de Brasilia, Brasília, Brazil. \\
   \fontsize{12}{12}\selectfont\textsuperscript{3}Departamento de Geometría y Topología, Universidad de Granada, Granada, Spain.
}
\date{}
\begin{document}
\maketitle


\begin{abstract}
    We consider weighted minimal surfaces in $\mathbb{H}^2 \times \mathbb{R}$, namely, surfaces that arise as critical points of a weighted area functional, where the corresponding weight depends only on the height function associated with the second factor $\mathbb{R}$. We refer to these surfaces as $\varphi$-minimal surfaces. For a class of increasing and concave weight potential functions $\varphi$, we focus our attention on $\varphi$-minimal surfaces that are invariant under isometries of $\mathbb{H}^2 \times \mathbb{R}$ that act trivially on $\mathbb{R}$ and act by hyperbolic translations in $\mathbb{H}^2.$ Through a detailed analysis of the corresponding three-dimensional system of ordinary differential equations, we present families of entire graphs, graphs over half-spaces of $\mathbb{H}^2$ that are asymptotic to vertical planes, winglike surfaces, and surfaces whose height function has no minimum. These families show that the hyperbolic-invariant setting is considerably richer than the rotational and parabolic settings. For the class of potential functions $\varphi$ we consider, they also allow us to classify all weighted minimal surfaces in $\mathbb{H}^2 \times \mathbb{R}$ that are invariant under hyperbolic translations. \\ \\
\end{abstract}

\noindent
\begin{it}2020 Mathematics Subject Classification.\end{it} 53A10; 53C42.

\medskip

\noindent
\begin{it}Keywords.\end{it} weighted minimal surfaces, product spaces, invariant surfaces, hyperbolic translations.

\section{Introduction}

Let $\varphi: \Omega \subseteq \mathbb{H}^2 \times \mathbb{R} \to \mathbb{R}$ be a differentiable function defined on an open subset $\Omega \subseteq \mathbb{H}^2 \times \mathbb{R}$. In this paper, we consider surfaces $\Sigma \subset \Omega$ that are critical points of the weighted area functional
$$
            A(\Sigma) = \int_\Sigma \text{e}^\varphi d\mu, 
$$
for every compactly supported variation of $\Sigma$. Equivalently, such surfaces are solutions of the prescribed mean curvature problem
\begin{equation}
\mathbf{H} = (\nabla \varphi)^\perp, \label{eq:main} 
\end{equation}
where $\mathbf{H}$ is the mean curvature vector of $\Sigma$ and $(\nabla \varphi)^\perp$ denotes the normal projection of $\nabla \varphi$ onto a unit normal $\eta$ of $\Sigma$. We call such surfaces $\varphi$-minimal surfaces. Besides being critical points of a weighted area functional, this terminology is also motivated by the fact that $\varphi$-minimal surfaces are also minimal surfaces in a Riemannian manifold endowed with a conformally changed metric $e^\varphi g$, where $g$ is the standard product metric of $\mathbb{H}^2 \times \mathbb{R}$ (see \cite{freitas2026geometry} for these equivalent properties).

Our primary interest lies in the case where the function $\varphi$ depends only on the height function with respect to the factor $\mathbb{R}$. Specifically, let $\sigma:\mathbb{H}^2\times \mathbb{R}\to\mathbb{R}$ be the height function $\sigma(p,z)=z$, and let $\partial_z=\nabla\sigma$ be the corresponding vertical vector field. For a smooth function $\varphi(z)$ depending only on $z$, it follows from equation \eqref{eq:main} that a surface $\Sigma\subset\mathbb{H}^2\times \mathbb{R}$ is a $\varphi$-minimal if its mean curvature $H$ satisfies
\begin{equation}\label{eq:phiminimal}
H = \dot{\varphi}(z)\,\langle \eta,\partial_z\rangle,
\end{equation}
where $\dot{\varphi}$ denotes the derivative of $\varphi$ with respect to $z$. Examples of $\varphi$-minimal surfaces include the classical minimal surfaces (when $\varphi(z)=c$, a constant), and the vertical translating solitons for the mean curvature flow (when $\varphi(z)=z$). Such solitons have recently been the subject of intense study, as we can see in the works \cite{bueno2018translating, lima-pipoli, de2019translating, bueno2024classgrimreapersmathbbh2timesmathbbr, tinaglia-pipoli}. 

Another motivation for considering this class of $\varphi$-minimal surfaces, beyond the classical and important examples mentioned above, is their recent study in Euclidean space, where $\varphi$ depends only on height. Again, when $\varphi(z)=z$, these surfaces become translating solitons of the mean curvature flow, which have been extensively studied in recent years. Martínez and Triviño \cite{martinez2022equilibrium} provided interesting advances in the study of such a class of $\varphi$-minimal surfaces in $\mathbb{R}^3$, with important contributions for symmetric surfaces. In addition, the study of translating solitons for the mean curvature flow has motivated the search for global results and asymptotic properties for this class of $\varphi$-minimal surfaces. See, for example, \cite{martinez2022equilibrium, mean-convex, catenary-cylinders, half-space}. Other interesting examples of $\varphi-$minimal surfaces in Euclidean space within this class are the singular minimal surfaces \cite{lopez2018invariant, lopez2019dirichlet, dierkes2021symmetric, weierstrass-martinez-trivino}, which occur when $\varphi(z) = \log(z)$.

In \cite{freitas2026geometry} we started the study of symmetric $\varphi$-minimal surfaces, motivated by previous works on symmetric translating solitons in $\mathbb{H}^2 \times \mathbb{R}$ (\cite{bueno2018translating, de2019translating, lima-pipoli} and symmetric $\varphi$-minimal surfaces in $\mathbb{R}^3$ \cite{martinez2022equilibrium}. For that, we consider isometries of $\mathbb{H}^2 \times \mathbb{R}$ that act trivially on $\mathbb{R}$ and act by parabolic translations, rotations, and hyperbolic translations $\mathbb{H}^2.$ We have shown that $\varphi$-minimal surfaces invariant under parabolic translation share properties analogous to a class of parabolic solitons in $\mathbb{H}^2\times \mathbb{R}$, and that rotationally symmetric $\varphi$-minimal surfaces in $\mathbb{H}^2\times \mathbb{R}$ behave similarly to rotational solitons in $\mathbb{H}^2 \times \mathbb{R}$ and their Euclidean counterparts.

The hyperbolic-invariant setting is geometrically richer, requires more involved analysis than the rotational and parabolic settings, and deserves special attention. That is the goal of this paper. In addition to the families mentioned above that we found in the rotationally and parabolically invariant setting, i.e., bowls, winglike surfaces and surfaces without minimum point at height, we present graphs over half-spaces of $\mathbb{H}^2$ that are asymptotic to vertical planes. Such a family appears here for the first time in the $\varphi$-minimal setting, having previously been known only for translating solitons, as we can see in \cite{bueno2024classgrimreapersmathbbh2timesmathbbr} and \cite{lima-pipoli}. This setting, however, presents significant challenges. On the one hand, the methods developed in the works \cite{bueno2024classgrimreapersmathbbh2timesmathbbr} and \cite{lima-pipoli} do not apply to general $\varphi$-minimal surfaces. On the other hand, the analysis we used to construct bowls, winglike surfaces, and surfaces without minimum at height is not enough to describe these graphs asymptotic to vertical planes. We overcome these difficulties through a detailed phase space analysis, considering not only the corresponding system of ordinary differential equations that describes the profile curve of those surfaces, but also an auxiliary system, and the behaviour of bowls and winglike surfaces as we vary their initial conditions.

First, we classify all $\varphi$-minimal surfaces in $\mathbb{H}^2\times \mathbb{R}$ that are symmetric about the vertical axis. We call these hyperbolic bowls by analogy with the rotational case. Second, we prove the existence of non‑symmetric vertical graphs and of a surface with winglike properties. Third, we determine a family of $\varphi$-minimal surfaces that have no minimum point in their height function. Finally, we establish a relation between the families of non-symmetric vertical graphs and winglike type surfaces, leading to our principal result: the existence of a vertical graph defined on a half‑space of $\mathbb{H}^2\times \mathbb{R}$ that is asymptotic to a vertical plane. Therefore, for surfaces invariant under hyperbolic translations, we obtain the following families, up to ambient isometries, leading to the classification of such invariant surfaces in the hyperbolic setting: 

\begin{itemize}
\item hyperbolic bowls (Theorem \ref{bowl hiperbolico});
\item hyperbolic vertical graphs whose height function does not attain a minimum, including the vertical-plane solution as a particular case  (Theorem \ref{sem minimo teorema});
\item hyperbolic entire vertical graphs whose profile curves are not symmetric with respect to the vertical line $x=0$ (Theorem \ref{teo grafico n simetricos});
\item hyperbolic winglike surfaces (Theorem \ref{winglike hiper});
\item hyperbolic vertical graphs contained in a half-space and asymptotic to a vertical plane (Propositions \ref{existenceasymptotic} and \ref{existenceasyntotic2}).
\end{itemize}

The paper is organized as follows. In Section 2 we recall preliminary facts about $\varphi$-minimal surfaces and set up the system of ODEs that describes surfaces invariant under hyperbolic translation. Section 3 treats the symmetric case (hyperbolic bowl). Section 4 is dedicated to prove the existence of hyperbolic $\varphi$-minimal surfaces that does not attain a minimum point at its height. Section 5 is 
devoted to non‑symmetric vertical graphs and the winglike surface. In Section 6 we show how these are related, prove the existence of the half‑space graph and complete the classification.

\section{Preliminaries}

Let $a \in \mathbb{R}$ and suppose that $\varphi : (a, \infty) \to \mathbb{R}$ satisfies
\begin{itemize}
    \item[(H1)] $\dot{\varphi} > 0$ and  $\ddot{\varphi} < 0$.
    \item[(H2)] {$\lim_{z \to \infty}z^2\ddot{\varphi}(z) = L_\infty > -\infty$.}
    \item[(H3)] $\lim_{z \to a^+}(z-a)\dot{\varphi}(z) = L_a \in ]0, \infty]$.
\end{itemize}

{\begin{remark} \label{rmk:behavior}
It is not difficult to prove that
\begin{itemize}
\item from  \rm{(H2)}, $\lim_{z\to\infty} z \dot{\varphi}(z) = L_\infty$ and $\lim_{z\to\infty}  \dot{\varphi}(z) = 0$,  
\item from \rm{(H3)},  $\lim_{z \to a^+} \frac{1}{\dot{\varphi}(z)} = 0$ and $\lim_{z \to a^+} \varphi(z) = -\infty$.
\end{itemize}
\end{remark}}

These hypotheses capture the essential features of $\varphi(z)=\ln z$ (where $a=0$), for which the prescribed mean curvature blows up as the surface approaches the height $z=a$ and decays to zero as $z\to\infty$. The same hypotheses also cover other singular-like behaviors.

Our goal is to describe $\varphi-$minimal surfaces invariant by one-parameter group of hyperbolic translations. Equivalently, surfaces foliated by equidistant curves. For simplicity, we will call such surfaces hyperbolic $\varphi-$minimal surfaces. To construct these surfaces, we follow the approach given in \cite{freitas2026geometry}. We consider the Lorentzian model of $\mathbb{H}^2 \times \mathbb{R}$ and we start with a planar curve
\begin{equation}\label{gamma}
  {  \Gamma(t) }= (\cosh x(t), \sinh x(t), 0, z(t)),
\end{equation}
parametrized by the arc length, which implies,
\begin{equation}\label{pca}
    x'^2 + z'^2 = 1.
\end{equation}

After applying the one-parameter group of hyperbolic translations, we obtain a surface $\Sigma$  parametrized by:
\begin{equation}\label{par hyperbolic}
    \psi(u,t) = \left(\cosh u \cosh x,\sinh x, \sinh u\cosh x, z  \right).
\end{equation}

A normal field on $\Sigma$ is giving by 
$$\eta(u,t)= -\left(z'\cosh u \sinh x, z'\cosh x, z' \sinh u \sinh x, -x' \right)$$
and the mean curvature is 
$$H(u,t
) = z'\tanh x +x'z'' - z'x''.$$
The condition that $\Sigma$ is $\varphi$-minimal, together with the expression for the mean curvature above, implies that the generating curve $\gamma$ must satisfy the following system of ordinary differential equations:
\begin{equation}\label{ODE system problem}
      \left\{  \begin{array}{lll}
        x'  &=& \cos \theta \\
         z' & =& \sin  \theta \\
         \theta ' &=&\dot{\varphi}(z) x'- z' \tanh x 
        \end{array}\right.
    \end{equation}

By the existence and uniqueness Theorem, given an initial condition, there exists a unique solution to this system. 

{A combination of Lemma 4.1 and Proposition 4.2} of \cite{freitas2026geometry}, we emphasize the following results.

\begin{prop}\label{maximal interval}
    If $x$, $z$ and $\theta$ are solutions of \eqref{ODE system problem} with the initial condition
    \begin{equation*}
        x(0) = x_0 \geq 0, \quad z(0) = z_0, \quad \theta(0) = \theta_0,
    \end{equation*}
    for $0\leq\theta_0<\frac{\pi}{2}$. Then, the solution is defined for all positive time.
\end{prop}
    
{Geometrically, $|x|$ provides the hyperbolic distance from a point on the profile curve to the point $P=(0,0,1)$, and its sign indicates its position in one of the half-spaces determined by the totally geodesic plane in $\mathbb{H}^2 \times \mathbb{R}$ given the geodesic $\alpha(s)=(\cosh(s),\sinh(s),0)$, which intersects the geodesic $\beta(s)=(\cosh(s), 0, \sinh(s))$ orthogonally at $P$. In turn, the function $z$ gives the height with respect the factor $\mathbb{R}$ and the function $\theta$ is the angle that the unit tangent vector field of the generating curve makes with the horizontal direction. Thus, understanding the behavior of these functions under different initial conditions supplies the information we need to describe the hyperbolic $\varphi-$minimal surfaces.}

Observe that if $\Sigma$ is a vertical graph $z = u(x)$,  then by (\cite{de2020hypersurfaces}, Equation (23)), we have that $u$ satisfies
\begin{equation} \label{edo graph}
    u'' = \left(1+u'^2\right)\left(\dot{\varphi}(u) - (\tanh x ) u'\right).
\end{equation}

In case  $\Sigma$ is a horizontal graph  $x = y(z)$. Then, $y$ satisfies the equation
\begin{equation}\label{graph hori}
        y'' = \left(1+y'^2\right)\left(\tanh y - \dot{\varphi}(z)y'\right).
    \end{equation}

\section{ The symmetric case} \label{symmetric case}

In this section, we analyze the case where the surface is symmetric about the {totally geodesic plane over the geodesic $\alpha(s)=(\cosh(s), \sinh(s),0)$. As we will see, this case} corresponds to problem \eqref{ODE system problem} with the initial conditions
\begin{equation} \label{initial cond symmetric}
    x(0)=0,\quad \theta(0)=0.
\end{equation}
 Note that by Proposition \ref{maximal interval} the solution is defined for all positive time. 

Initially, define the right branch $\gamma^+$ and left branch $\gamma^-$ of $\gamma$, respectively, as $\gamma|_{(0, \infty)}$ and $\gamma|_{(T_{\min}, 0)}$, where $(T_{\min}, 0)$ is the maximal interval solution in a negative time. To prove the existence of a symmetric vertical graph it is necessary the following auxiliary Lemma.

\begin{lemma} \label{lema1}
    Let $\gamma = (x, z, \theta)$ be a solution of (\ref{ODE system problem}) with an initial condition 
    \begin{equation}\label{initial cond.}
        x(0) = x_0 \geq 0, \quad z(0) = z_0 \in (a, \infty), \quad \text{and} \quad \theta(0) = \theta_0,
    \end{equation}
    for $0\leq \theta_0 < \frac{\pi}{2}$. Then,  $x$ and $z$ increase for $t > 0$. Moreover,
    \begin{equation*}
        \lim_{t \to \infty} x = \infty, \quad \lim_{t \to \infty} z = \infty, \quad \lim_{t \to \infty}\theta = 0.
    \end{equation*}
\end{lemma}
\begin{proof}
    Observe that with this initial condition, the functions $x$ and $z$ increase near $t = 0$. Suppose, by contradiction, that there exists a $t_0 > 0$ such that $x'(t_0) = 0$, then it should be a maximum point or a saddle point. However,
    $$x''(t_0) = -\theta'(t_0) \sin \theta (t_0) = \sin^2\theta\left(t_0\right)\tanh(x) > 0.$$
    Therefore, $x' > 0$ for all $t > 0$. Similarly, $z' > 0$. Consequently, $0 < \theta(t) < \pi/2$, and, by Proposition \ref{maximal interval}, the solution is defined for all positive time. 

    Next, observe that the sign of $\theta'(0)$ cannot be determined if $\theta(0) > 0$. However, it is still possible to prove that $\theta$ has at most one critical point. Indeed, suppose, by contradiction, that $t_0$ and $t_1$ are critical points of $\theta$ with $t_0 < t_1$. Then, 
    \begin{equation} \label{second derivate}
        \theta '' (t_0) = \cos \theta \sin \theta \left(\ddot{\varphi}(z(t_0))-\frac{1}{\cosh^2 x(t_0)} \right) < 0.
    \end{equation}
 Without loss of generality, suppose that $t_1$ is the first critical point of $\theta$ after $t_0$. In this case, $\theta''(t_1) \geq 0$. On the other hand, analogously to \eqref{second derivate}, $\theta''(t_1) < 0$, which is a contradiction.
Moreover,
$$\theta'' = \ddot{\varphi}(z)\sin\theta\cos\theta
-\frac{\sin\theta\cos\theta}{\cosh^2(x)}
-\theta'\left(\dot{\varphi}(z)\sin\theta+\tanh(x)\cos\theta
\right).$$
Since and
$0<\theta<\pi/2$, we obtain
$$ |\theta''| \leq |\ddot{\varphi}(z)|+\frac{x'}{\cosh^2(x)}+\bigl(\dot{\varphi}(z_0)+1\bigr)|\theta'|.
$$
Since $\ddot{\varphi}<0$, we have
$$
\int_0^t |\ddot{\varphi}(z(s))|z'(s)\,ds
=
\int_{z_0}^{z(t)}|\ddot{\varphi}(u)|\,du
=
\dot{\varphi}(z_0)-\dot{\varphi}(z(t))
\leq\dot{\varphi}(z_0).
$$
Moreover,
$$ \int_0^t\frac{x'(s)}{\cosh^2(x(s))}\,ds
=
\tanh(x(t))-\tanh(x_0)
\leq1. $$
Since $\theta$ is bounded and has at most one critical point,
$$ \int_0^\infty |\theta'(s)|\,ds<\infty. $$

Therefore,
$$
\int_0^\infty|\theta''(s)|\,ds<\infty,
$$
and hence $\theta'$ is uniformly continuous.

Furthermore, since $\theta$ is bounded and eventually monotone,
$\lim_{t\to\infty}\theta(t)$ exists, and therefore
$$
\int_0^\infty\theta'(s)\,ds
=
\lim_{t\to\infty}\theta(t)-\theta(0)
$$
is finite. Thus, by Barbalat's lemma (see
\cite{khalil2002nonlinear}),
\begin{equation}\label{limi theta}
\lim_{t\to\infty}\theta'(t)=0.
\end{equation}

   Before analyzing the limit  $\lim_{t \to \infty} \theta(t)$, we first investigate $\lim_{t \to \infty} x(t)$ and $\lim_{t \to \infty} z(t)$. So, suppose two different cases: $\theta' > 0$ for all $t > 0$ or there exists $t_0 > 0$ such that $\theta'(t) < 0$ for all $t > t_0$. In the first case, $z'' > 0$ implying that $z$ is an unbounded function. Suppose, by contradiction that $\lim_{t \to \infty} x(t) = x^+ $ is bounded, then $\lim_{t \to \infty} x'(t) = 0$. Consequently, 
   $$\lim_{t \to \infty}\theta'(t) = -\tanh x^+<0.$$
   This contradicts \eqref{limi theta}. 

   In the second case, $x'' > 0$ and it follows that $\lim_{t \to \infty} x(t) = \infty$. Suppose, by contradiction, that $z$ is bounded from above and $$\lim_{t \to \infty} z(t) = z_1,$$ by the ODE system (\ref{ODE system problem}), we have that
    \begin{align*}
        \theta ' & > x'\dot{\varphi}(z)-z'.
    \end{align*}
    Integrating this inequality, we obtain that
    \begin{align*}
        \theta + z - z_0 &> \int_0^t(x'\dot{\varphi}(z)) \\
            &=  x\dot{\varphi}(z) - \int_0^t xz'\overset{..}{\varphi}(z) \\
            &> x\dot{\varphi}(z).
    \end{align*}
    Then
    \begin{equation} \label{relation x y z}
        \frac{z +c_1}{\dot{\varphi}(z)} > x,
    \end{equation}
    where $c_1$ is a constant. Since $\lim_{t \to \infty} x(t) = \infty$, we have a contradiction. 

    In conclusion, $x$ and $z$ are unbounded functions. Then, by \eqref{limi theta},
    $$ 0 = \lim_{t \to \infty}\dot{\varphi}(z) = \lim_{t \to \infty}\tan \theta,$$
    implying that $\theta$ tends to zero.
\end{proof}
 {From the proof of Lemma \ref{lema1}, we have the following consequence:}
\begin{cor}\label{corol1}
If $\theta_0 = 0$, then $\theta$ has a unique critical point in $(0, \infty)$, which is a maximum point.
\end{cor}

By using Lemma \ref{lema1}, we can prove the existence of a surface, which we will call the \textit{ $\varphi$-minimal hyperbolic bowl}, motivated by \cite{de2019translating}, since this surface has properties similar to the translating soliton Bowl (see \cite{bueno2018translating}).  
{\begin{remark}\label{Pi/2} If $x_0>0$, Lemma \ref{lema1} also holds for $\theta_0=\frac{\pi}{2}$.
\end{remark}}
\begin{theorem}\label{bowl hiperbolico}
    There exists a $\varphi$-minimal surface in $\mathbb{H}^2\times \mathbb{R}$, properly embedded, foliated by the equidistant curves, which is {an entire vertical } graph, symmetric about the vertical axis, intersecting this axis orthogonally,  {with  unbounded height function and with angle function tending to $0$ at infinity.}
\end{theorem}
\begin{proof}
    Consider the problem (\ref{ODE system problem}) with the initial conditions 
    \begin{equation} \label{simetrica}
        x(0) = 0, \quad z(0) = z_0, \quad \theta(0) = 0.
    \end{equation}
    Suppose that the maximal solution is defined on the interval $(t_{\min}, t_{\max})$, and observe that it is symmetric about the vertical axis. Indeed, for $t \in (t_{\min}, 0)$, the functions 
    $\tilde{x}(t) = -x(-t)$, $\tilde{z}(t) = z(-t)$ and $\tilde{\theta}(t) = -\theta(-t)$ are solutions of (\ref{ODE system problem})-(\ref{simetrica}). Therefore, it is sufficient to describe the solution's behavior for positive time. By Lemma \ref{lema1}, we conclude that $x'> 0$, which implies the surface is a vertical graph. Moreover, we conclude from the limits that
    { $$\lim_{t \to \infty }z(t) = \infty, \quad  \quad \lim_{t \to \infty}x(t) = \infty \quad \text{and}\quad  \lim_{t \to \infty}\theta(t) = 0.$$ In particular,}  the height function is unbounded and the surface is a complete graph  {whose angle function tends to $0$ at infinity.} 
\end{proof}

\begin{figure}[ht]
  \centering
  \begin{minipage}[b]{0.45\textwidth}
    \includegraphics[scale=0.75]{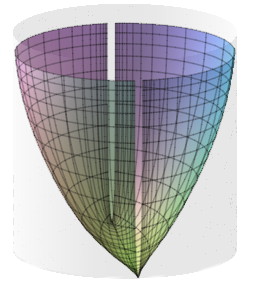}
    \caption{$\varphi$-minimal hyperbolic bowl for $\varphi(z) = \text{e}^{-z}$}
  \end{minipage}
  \hfill
  \begin{minipage}[b]{0.43\textwidth}
    \includegraphics[scale=0.75]{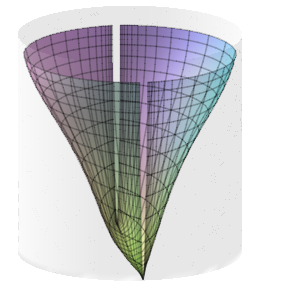}
    \caption{$\varphi$-minimal hyperbolic bowl for $\varphi(z) = \frac{1}{z}$}
  \end{minipage}
\end{figure}

Another example of $\varphi-$minimal surfaces is the vertical plane, which is obtained when the initial condition is $x(0) = 0$ and $\theta(0) = \pm \frac{\pi}{2}$, described in the following proposition.

\begin{prop}
    If $\gamma = (x, z,\theta)$ is a solution of (\ref{ODE system problem}) with the initial conditions
    \begin{equation*}
        x(0) = 0, \quad z(0) = z_0 \quad \text{and} \quad \theta(0) = \frac{\pi}{2} \text{( or $-\frac{\pi}{2}$)}.
    \end{equation*}
    Then, for all $t \in \mathbb{R}$,
    $x(t) = 0$,  $\theta(t) = \frac{\pi}{2} \left( \text{or} -\frac{\pi}{2}\right)$ and $z(t) = z_0+t$ (or $z(t)=z_0-t$).
\end{prop}

\section{The \texorpdfstring{$\varphi$}{phi}-minimal surface without minimum}

This section is dedicated to prove existence of a $\varphi$-minimal
surface that does not attain a minimum point at its height. This
result also addresses another relevant question regarding the
solutions to \eqref{ODE system problem}, specifically their
asymptotic behavior near the singular point $z=a$.

Assume that $x=y(z)$ a solution of \eqref{graph hori} in
$]a,a+\delta[$ that extends to $C^1([a,a+\delta[)$ with
$y(a)=x_0$ and $y'(a)=y_0$. To establish the existence
of solutions of this type, we will follow the approach established
in \cite{lopez2018invariant} defining an operator on
$C^1([a,a+\delta])$ and showing that it admits a fixed point
which corresponds to a solution of \eqref{graph hori} satisfying
the initial conditions $y(a)=x_0$ and $y'(a)=y_0$. Indeed,
we observe that they are solutions to the following degenerate
initial value problem:
\begin{equation}\label{sis equ}
\begin{cases}
\displaystyle
\left(
\frac{y'}{\sqrt{1+y'^2}}e^{\varphi(z)}
\right)'
=
\frac{\tanh y}{\sqrt{1+y'^2}}e^{\varphi(z)},
\\[1.2em]
y(a)=x_0,\quad y'(a)=y_0.
\end{cases}
\end{equation}

Then, defining the functions
$f:\mathbb{R}_a^+\times\mathbb{R}\to\mathbb{R}$
and $g:\mathbb{R}\to\mathbb{R}$ by
$$
f(s,t)=\frac{\tanh s}{\sqrt{1+t^2}}
\qquad\text{and}\qquad
g(t)=\frac{t}{\sqrt{1+t^2}},
$$
the problem \eqref{sis equ} can be rewritten as
\begin{equation}\label{equiv ode}
\begin{cases}
\left(g(y')e^\varphi\right)'=e^\varphi f(y,y'),
\\[0.8em]
y(a)=x_0,\quad y'(a)=y_0.
\end{cases}
\end{equation}

\begin{lemma}\label{initial derivative zero}

 Let $y=y(z)$, $z\in]a,a+\delta[$ be  a solution of \eqref{equiv ode} and suppose that there exist $\lim_{z\to a}y(z) = x_0$ and  $\lim_{z\to a} y'(z)=L_0\in[-\infty,\infty]$. Then $L_0=0$.
\end{lemma}

\begin{proof}
From \eqref{equiv ode}, for any $\epsilon>a$, we have
$$
g(y'(z))e^{\varphi(z)}
-
g(y'(\epsilon))e^{\varphi(\epsilon)}
=
\int_\epsilon^z
f(y(t),y'(t))e^{\varphi(t)}\,dt.
$$
Letting $\epsilon\to a$ in the previous expression and taking
into account Remark \ref{rmk:behavior} and the fact that $g$ is bounded, we obtain
$$
g(y'(z))
=
e^{-\varphi(z)}
\int_a^z f(y(t),y'(t))e^{\varphi(t)}\,dt.
$$
But then, by L'Hôpital's rule, we have
$$
0\leq |g(y_0)|
=
\lim_{z\to a}|g(y'(z))|
\leq
\lim_{z\to a}
\frac{\displaystyle\int_a^z e^{\varphi(t)}\,dt}
     {e^{\varphi(z)}}
=
\lim_{z\to a}\frac{1}{\dot{\varphi}(z)}
=
0,
$$
which concludes the proof.
\end{proof}

Consider the operator
$I:C^1([a,\infty))\to C^1([a,\infty))$ defined by
$$
I(y)(s)
=
e^{-\varphi(s)}
\left(
\int_a^s f(y,y')e^{\varphi(t)}\,dt
\right).
$$
Since, from Remark 2.1,
$\lim_{s\to a}\varphi(s)=-\infty$, then by L'Hôpital's rule
we have,
\begin{equation}\label{limit integral operator}
\lim_{s\to a}
\frac{\displaystyle\int_a^s e^{\varphi(t)}\,dt}
     {e^{\varphi(s)}}
=
\lim_{s\to a}
\frac{e^{\varphi(s)}}
     {\dot{\varphi}(s)e^{\varphi(s)}}
=
0,
\end{equation}
and there exists $\delta>0$ such that $|I|<1$ for all
$s\in(a,a+\delta)$. Consequently, the operator
$S:C^1[a,a+\delta]\to C^1[a,a+\delta]$ defined as
$$
S(y(z))
=
\int_a^z
g^{-1}\left(
e^{-\varphi(s)}
\int_a^s f(y,y')e^{\varphi(t)}\,dt
\right)\,ds
+x_0
$$
is well defined.

\begin{prop}\label{contraction}
The operator $S$ is a contraction.
\end{prop}

\begin{proof}
{Since the functions $f$ and $g^{-1}$ are Lipschitz continuous
in the domains
$[x_0-\epsilon,x_0+\epsilon]\times[-\epsilon,\epsilon]$
and $[-\epsilon,\epsilon]$, respectively, provided $\epsilon$,
$0<\epsilon<1$,  then for all $y,u\in B(x_0,\epsilon)$}
and for $z\in[a,a+\delta]$, we have the estimate
\begin{align*}
|Sy(z)-Su(z)|
&\leq
c_1\int_a^z e^{-\varphi(s)}
\left|\int_a^s e^{\varphi(t)}
\left|f(y,y')-f(u,u')\right|\,dt\right|\,ds
\\
&\leq
c_2\left(
\|y-u\|_\infty+\|y'-u'\|_\infty
\right)
\int_a^z\int_a^s
e^{-\varphi(s)}e^{\varphi(t)}\,dt\,ds,
\end{align*}
where $c_1$ and $c_2$ are positive constants.

From \eqref{limit integral operator}, we may choose $\delta>0$
small enough such that the integral
$$
\int_a^s e^{-\varphi(s)}e^{\varphi(t)}\,dt
$$
is bounded. Hence, there exists a positive constant $c_3$
such that
$$
|Sy(z)-Su(z)|
\leq
c_3\left(
\|y-u\|_\infty+\|y'-u'\|_\infty
\right),
$$
and we conclude that for $\delta>0$ small enough, $S$ is a
contraction in the closed ball {$\overline{B(x_0,\epsilon)}$}
in $C^1([a,a+\delta])$.
\end{proof}

\begin{theorem}\label{existence}
There exists a unique solution $y\in C^2([a,a+\delta])$ of
the integral equation
$$
y(z)
=
\int_a^z
g^{-1}\left(
e^{-\varphi(s)}
\left(
\int_a^s f(y(t),y'(t))e^{\varphi(t)}\,dt
\right)
\right)\,ds
+x_0,
$$
which is a solution of
\begin{equation}\label{equiv ode zero}
\begin{cases}
\left(g(y')e^\varphi\right)'=e^\varphi f(y,y'),
\\[0.8em]
y(a)=x_0,\quad y'(a)=0.
\end{cases}
\end{equation}
\end{theorem}

\begin{proof}
Proposition \ref{contraction} ensures the existence and
uniqueness of the solution $y(z)$ in $C^1([a,a+\delta])$.

Moreover, since
$$
g(y')
=
e^{-\varphi(z)}
\int_a^z f(y,y')e^{\varphi(t)}\,dt,
$$
we have by L'Hôpital's rule,
$$
\begin{aligned}
\lim_{z\to a^+}\frac{g(y')}{z-a}
&=
\lim_{z\to a^+}
\frac{\displaystyle\int_a^z f(y,y')e^{\varphi(t)}\,dt}
     {e^{\varphi(z)}(z-a)}
\\
&=
\lim_{z\to a^+}
\frac{f(y,y')}{\dot{\varphi}(z)(z-a)+1}
=
\frac{\tanh(x_0)}{1+L}.
\end{aligned}
$$
On the other hand, from definition of $g$,
$$
\frac{\tanh(x_0)}{1+L}
=
\lim_{z\to a^+}\frac{g(y')}{z-a}
=
{\lim_{z\to a^+}\frac{g(y')\sqrt{1+y'^2}}{z-a}}
=
\lim_{z\to a^+}\frac{y'(z)}{z-a}
=
y''(a),
$$
and $y'$ is differentiable in $a$. The continuity of $y''$
in $a$ follows by a straight forward computation letting
$z\to a$ in \eqref{graph hori}.
\end{proof}

\begin{remark}
The previous theorem is not true without hypothesis (H3).
Indeed, taking
$$
\dot{\varphi}(z)
=
\frac{1}{z(3+\sin(\log(z))},
\qquad z\in]0,\infty[,
$$
Proposition \ref{contraction} yields the existence of a
solution to \eqref{equiv ode zero} that can not be extended
to $C^2([0,\delta])$.
\end{remark}

{We are now in position to establish the theorem that presents the family of $\varphi$-minimal surfaces whose height function has no minimum.}

\begin{theorem}\label{sem minimo teorema}
For any $x_0$, there exists a hyperbolic $\varphi$-minimal
surface $\Sigma$, regular up to the boundary$\partial\Sigma$,
contained in the half-space $\{z>a\}$ of $\mathbb{H}^2\times \mathbb{R}$, satisfying
that $(\cosh(x_0),\sinh(x_0),0,a)\in\partial\Sigma$ and for
which the height function and does not attain a minimum on
$\Sigma$. Moreover, both the distance and the height function
are unbounded functions.
\end{theorem}

\begin{proof}
If $x_0=0$, from \eqref{ODE system problem} and \eqref{equiv ode zero}, it is clear
that the half-vertical plane $\{x=0,z>a\}$ gives the example
we are looking for.

If $x_0\neq 0$, then by Theorem \ref{existence}, since the
fixed point exists for arbitrarily small $\delta$, we can
extend the solution $y(z)$ to \eqref{equiv ode zero} to a
maximal interval of existence. Consequently, there exists a
hyperbolic $\varphi$-minimal surface that does not attain
a minimum at its height because $y'$ is strictly monotonic
when $z\to a$.

To analyze the global behavior, rewrite $y$ as
$y(z)=x(z(t))$ where $x$ and $z$ satisfy the first-order
system \eqref{ODE system problem}. In a small neighborhood
of $t=0$, we may choose initial conditions
$$
x(0)=x_0,\qquad z(0)=z_0,\qquad \theta(0)=\theta_0,
$$
such that $(x_0,z_0,\theta_0)$ is arbitrarily close to
$(x_0,a,\pi/2)$. The correspondence is given by
$y(z_0)=x_0$ and $\tan\theta_0=y'(z_0)$; the limit
$y'(a)=0$ implies $\theta_0\to\pi/2$.

Consider $x_0>0$, we claim that $\theta_0<\frac{\pi}{2}$.
Applying Lemma \ref{lema1} will then yield the desired result.

Suppose, for the sake of contradiction, that
$\theta_0>\frac{\pi}{2}$. In this case, $\theta$ increase
near $t=0$ and $\theta'(0)\geq 0$. However, in a small
neighborhood of $t=0$, system \eqref{ODE system problem}
yields
$$
0>\dot{\varphi}(z)\cos\theta
\geq \sin\theta\tanh x>0,
$$
which is a contradiction. Thus, if $x_0>0$, then
$\theta_0<\frac{\pi}{2}$. The case where $x_0<0$ follows
by the symmetry
$$
\widetilde{x}=-x,\qquad
\widetilde{z}=z,\qquad
\widetilde{\theta}=\pi-\theta.
$$
This completes the proof.
\end{proof}

\begin{figure}[ht]
    \centering
    \begin{minipage}[b]{0.45\textwidth}
        \includegraphics[scale=0.4]{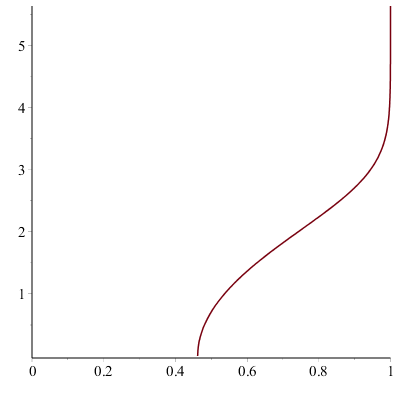}
    \end{minipage}
    \hfill
    \begin{minipage}[b]{0.43\textwidth}
        \includegraphics[scale=0.5]{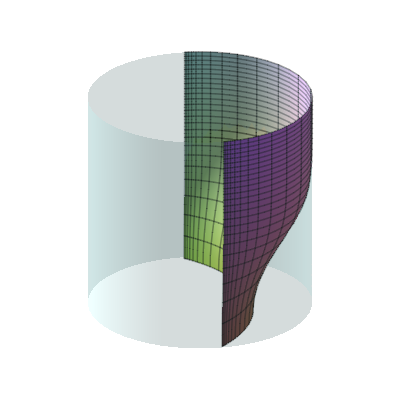}
    \end{minipage}
    \caption{$\varphi$-minimal hyperbolic without minimum
    for $\dot{\varphi}(z)=\frac{1}{z}$.}
    \label{fig:sem-minimo}
\end{figure}

Figure \ref{fig:sem-minimo} illustrates this family of
$\varphi$-minimal surfaces in $\mathbb{H}^2\times\mathbb{R}$
when $\varphi(z)=\ln z$.

To study which kind of solutions of \eqref{ODE system problem} generate the examples described in Theorem \ref{sem minimo teorema}, we first prove the following 
\begin{lemma}\label{lema3}
    Let $x$, $z$ and $\theta$ be a solution of (\ref{ODE system problem}) with the initial condition
    \begin{equation}\label{initial cond. non symmetric}
    x(0) = x_0 \geq 0, \quad z(0) = z_0 > a, \quad\theta(0) = -\theta_0,
    \end{equation}
    where $0<\theta_0<\frac{\pi}{2}$.
    Then, the solution is defined for all positive time, $x$  is strictly increasing, and  $z$ has a unique minimum point. Moreover,{$$\lim_{t \to \infty}x(t) = \infty, \quad \lim_{t \to \infty}z(t) = \infty \quad \text{and} \quad \lim_{t \to \infty} \theta(t) = 0$$}
\end{lemma}
\begin{proof}

Denote by $(0, T_{\max})$ the maximal interval of existence of the solution to \eqref{ODE system problem}–\eqref{initial cond. non symmetric}. Then, $x'(0) > 0$, $z'(0) < 0$ and, by an argument analogous to that used in Lemma \ref{lema1}, one can show that $x'(t) > 0$ for all $t > 0$.

If $z'(t) < 0$ for all $t \in (0, T_{\max})$, then, $-\frac{\pi}{2} < \theta(t) < 0$ and we claim that 
 $\lim_{t \to T_{\max}} z(t) = a$.

Otherwise, $\lim_{t \to T_{\max}} z(t) > a$, $\lim_{t \to T_{\max}} x(t) = \infty$, the maximal interval is unbounded and, from Barbalat’s Lemma, it follows that
$$\lim_{t \to T_{\max}} z'(t) = 0 \quad \text{and} \quad \lim_{t \to T_{\max}} x'(t) = 0,$$
which is a contradiction.  

Since  $\theta' = \cos \theta \dot{\varphi}(z) - \sin \theta \tanh x > 0$, there exists $\theta_L=\lim_{t\to T_{max}} \theta(t)\in (-\frac{\pi}{2},0]$ and then, from Lemma \ref{initial derivative zero}, and having in mind that $(x(t), z(t))$ is a horizontal graph near $z = a$,  we also have that $\lim_{t \to T_{\max}} x'(t) = 0$, which is also a contradiction. 

Consequently, there exists a first instant $t_0 > 0$ such that $z'(t_0) = 0$ and since $z''(t_0) = \cos^2 \theta(t_0) \dot{\varphi}(z(t_0)) > 0$, we conclude that $z(t_0)$ is a minimum point. If $z$ has another critical point, then there exists $t_1 > t_0$ such that $z'(t_1) = 0$ and $z''(t_1) \leq 0$. But by a direct computation, $z''(t_1) > 0.$ Thus, $z$ has a unique critical point. As proved in Lemma 4.1 in \cite{freitas2026geometry}, $z$ is defined for all positive time and the limits are verified as in Lemma \ref{lema1}.
\end{proof}
\begin{theorem}\label{unicas}
Let $\gamma=(x,z,\theta)$ be a solution of
\eqref{ODE system problem} on a maximal interval
$I_M=(T_{\min},T_{\max})$ satisfying that there exists $t_0\in I_M$ with 
\begin{equation} \label{conditionminimum} x(t_0)=x_0<0, \qquad -\frac{\pi}{2}<\theta<0, \quad \text{on} \quad (t_0,T_{max}).\end{equation}
Then, the hyperbolic $\varphi$-minimal surface generated
by $\gamma$ is one of the examples described in
Theorem \ref{sem minimo teorema}.
\end{theorem}
\begin{proof}
From \eqref{conditionminimum}, we have that $x'>0$ and $z'<0$ on $(t_0,T_{\max})$. Thus, there exist $x_L=\lim_{t\to T_{\max}}x(t)$ and $z_L=\lim_{t\to T_{\max}}z(t)$. Moreover, if $\theta'(t_1)=0$ for some $t_1\geq t_0$, then 
$$\theta''(t_1) =\ddot{\varphi} \cos \theta \sin \theta - \frac{1}{\cosh^2 x} \cos \theta \sin \theta >0,$$
and $\theta$ has at most one critical point on $(t_0,T_{\max})$. Consequently, $\theta$ is strictly monotonic near $T_{\max}$ and  there exists $\theta_L=\lim_{t\to T_{\max}}\theta(t)\in [-\frac{\pi}{2},0]$.

From Lemma \ref{lema3}, $x_L\leq 0$ and since $I_M$ is maximal, we also deduce that $z_L=a$. In particular, near $T_{\max}$, $(x,z)$ is a horizontal graph $x=y(z)$ with $x(a)=x_L$ and $y'(a)=\cot \theta_L$. The proof concludes from Lemma \ref{initial derivative zero} and Theorem \ref{existence}.
\end{proof}

\section{Non-symmetric case and winglike surfaces}

This section deals with the existence of two different classes of surfaces: one consisting of winglike surfaces and the other composed of vertical graphs not necessarily symmetric with respect to the vertical axis.

These surfaces will be obtained by analyzing the solutions of \eqref{ODE system problem} with initial conditions \eqref{initial cond.}. The right branch $(t > 0)$, of the solution $\gamma = (x,z,\theta)$ of \eqref{ODE system problem}-\eqref{initial cond.} has already been described in Lemma \ref{lema1}. For an analysis of the left branch $(t < 0)$, we will introduce suitable changes of variables depending on the case we wish to describe. Both vertical graphs and winglike surfaces verify that their height function attains a minimum. 

\begin{figure}[ht]
    \centering
    \includegraphics[scale=0.75]{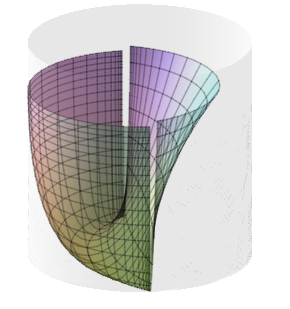}
    \caption{$\varphi$-minimal graph for $\varphi(z) = -e^{1/z}$}
    \label{fig:enter-label}
\end{figure}

Combining Lemma \ref{lema1} and Lemma \ref{lema3}, it is possible to describe a non-symmetric vertical graph, shown in Figure \ref{fig:enter-label}, in the following theorem.

\begin{theorem} \label{teo grafico n simetricos}
    There exists a $\varphi$-minimal surface, properly embedded, foliated by equidistant curves, that is an entire vertical graph, non-symmetric about the vertical axis. The height function is unbounded   { and has  a minimum point, while the angle function tends to $0$ at infinity}.  
\end{theorem}
\begin{proof}
    Consider $x$, $z$ and $\theta$ solutions of \eqref{ODE system problem} with the initial condition 
    $$x(0) = 0, \quad z(0) = z_0 \quad \theta(0) = \theta_0,$$
    where $0<\theta_0<\frac{\pi}{2}$.  By Lemma \ref{lema1}, the right branch of the profile curve 
    $$\gamma(t) = (\cosh x(t), \sinh x(t), 0, z(t)),$$
    describes a vertical graph with   {unbounded height and distance functions and with an angle function tending to $0$ at infinity}.
    For the left branch of $\gamma$, consider the variable change
        \begin{equation*}
    \tilde{x}(t)= -x(-t), \quad \tilde{z}(t) = z(-t), \quad \text{and} \quad \tilde{\theta}(t) = -\theta(-t).
\end{equation*}
Note that $\tilde{x}$, $\tilde{z}$ and $\tilde{\theta}$ are solutions of (\ref{ODE system problem}) with the initial conditions
\begin{equation*}
    \tilde{x}(0)= 0, \quad \tilde{z}(0) = z_0, \quad \text{and} \quad \tilde{\theta}(0) = -\theta_0.
\end{equation*}
    Then, by Lemma \ref{lema3}, the left branch is a vertical graph where the height function has a unique minimum point, with   {unbounded height and distance functions and with an angle function going to $0$ at infinity}.
\end{proof}

Another $\varphi-$minimal surface, with different properties, is obtained when the solution $\gamma$ of (\ref{ODE system problem}) is chosen with the initial conditions  $\theta(0) = \frac{\pi}{2}$ and $x(0)\neq 0$. The most interesting property about the solution  $\gamma$ is that each branch is a vertical graph, but it is not a vertical graph at $\gamma(0)$.
In this case, the surface is described as winglike, similarly to \cite{bueno2018translating} and \cite{martin2015topology}. 

\begin{figure}[ht]
  \centering
  \begin{minipage}[b]{0.45\textwidth}
    \includegraphics[scale=0.5]{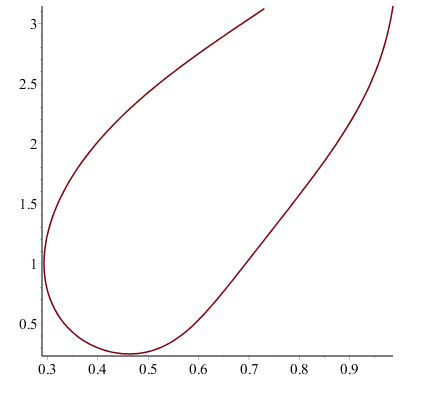}
  \end{minipage}
  \hfill
  \begin{minipage}[b]{0.43\textwidth}
    \includegraphics[scale=0.75]{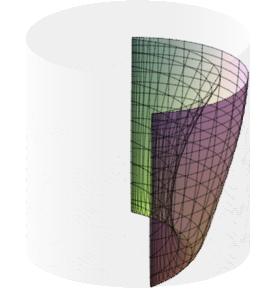}
  \end{minipage}
  \caption{hyperbolic $\varphi$-minimal winglike, for $\dot{\varphi}(z) = \frac{1}{z}$}\label{fig winglike}
\end{figure}

\begin{lemma}\label{cor-winglike-case}
    Let $x$, $z$ and $\theta$ be a solution of (\ref{ODE system problem}) with the initial condition
    \begin{equation*}
    x(0) = x_0 > 0, \quad z(0) = z_0 > a, \quad\theta(0) = -\frac{\pi}{2}. 
    \end{equation*}
    Then, the solution is defined for all positive time, $x$ is strictly increasing, and and $z$ has a unique minimum point. Moreover,   {$$\lim_{t \to \infty}x(t) = \infty, \quad \lim_{t \to \infty}z(t) = \infty \quad \text{and} \quad \lim_{t \to \infty} \theta(t) = 0$$}
\end{lemma}
\begin{proof}
    Observe that $x'(0) = 0$ and $x''(0) > 0$. Then, $x$ increases near $t>0$ and the proof follows as in Lemma \ref{lema3}.
\end{proof}
\begin{theorem} \label{winglike hiper}
There exists a $\varphi$-minimal surface, properly embedded, foliated by equidistant curves which is a non-convex, vertical bigraph.  This surface   {is contained in a half-space of $\mathbb{H}^2\times \mathbb{R}$, the height and distance functions are unbounded above but they have both  a minimum and the angle function varies from $-\pi$ to $0$ at infinity. }.
\end{theorem}
\begin{proof}
    Consider a solution of \eqref{ODE system problem} with the initial condition
    \begin{equation*}
        x(0) = x_0 > 0, \quad z(0) = z_0, \quad \theta(0) = \frac{\pi}{2}. 
    \end{equation*}
From Remnark \ref{Pi/2}, the right branch of the profile curve 
    $$    {  \Gamma(t) } = (\cosh x(t), \sinh x(t), 0, z(t)),$$
    describes a vertical graph with unbounded height and distance,   { and  with angle function tending to $0$ at infinity}

    For $t < 0$, consider a variable change 
    \begin{equation*}
        \tilde{x}(t) = x(-t), \quad \tilde{z}(t) = z(-t), \quad \tilde{\theta}(t) = \theta(-t) - \pi.
    \end{equation*}
    Note that $\tilde{x}$, $\tilde{z}$ and $\tilde{\theta}$ are solutions of \eqref{ODE system problem} with the initial condition 
    \begin{equation*}
        \tilde{x}(0) = x_0, \quad \tilde{z}(0) = z_0, \quad \tilde{\theta}(0) = -\frac{\pi}{2}.
\end{equation*}
Therefore, by Lemma \ref{cor-winglike-case}, $\tilde{z}$ has a minimum point and the left branch of $\gamma$ describes a vertical graph, with unbounded height and distance,   { and  with angle function tending to $-\pi$ at infinity} .   
\end{proof}
The figure \ref{fig winglike} shows the winglike surface and its profile curve for $\varphi(z) = \ln z$.


\section{Half-space graphs and classification}
In this section we complete the list of examples and classify all $\varphi$-minimal surfaces invariant under hyperbolic translations.
\subsection{Solutions with minimum}
First, we investigate, from a geometric perspective, the trajectories of \eqref{ODE system problem} that possess a minimum with respect to the height function.
\begin{theorem}\label{withminimum} Let $\gamma=(x,z,\theta)$ be a solution of \eqref{ODE system problem} with $\gamma(0)=(x_0,z_0,0)$. Then the $\varphi$-minimal surface $\Sigma$ generated by $\gamma$ is one of the following examples,
\begin{itemize}
\item[(A)] A symmetric  entire vertical $\varphi$-graph.
\item[(B)] A non-symmetric entire vertical $\varphi$-graph.
\item[(C)] A winglike $\varphi$-minimal surface.
\item[(D)] A vertical $\varphi$-graph in the half-space $\{x>0\}$ (or $\{x<0\})$ asymptotic to he vertical plane $\{x=0\}$ and with unbounded height function  at infinity.
\end{itemize}
\end{theorem}
\begin{proof}
If $x_0=0$, then, from Theorem \ref{bowl hiperbolico}, $\Sigma$ is a symmetric entire vertical graph.
\\
If $x_0\neq 0$  , we may assume without loss of generality that  $x_0>0$; otherwise, it suffices to consider the solution $\hat{\gamma} =(\hat{x} ,\hat{z} ,\hat{\theta} )$ given by 
    \begin{equation*}
    \hat{x}(t)= -x(-t), \quad \hat{z}(t) = z(-t), \quad \text{and} \quad \hat{\theta}(t) = -\theta(-t).
\end{equation*}
Then, from Lemma \ref{lema1}, $x$ and $z$ are strictly increasing for $t > 0$. Moreover,
    \begin{equation*}
        \lim_{t \to \infty} x = \infty, \quad \lim_{t \to \infty} z = \infty, \quad \lim_{t \to \infty}\theta = 0.
    \end{equation*}
To study the left branch of $\gamma$ we set 
    \begin{equation*}
    \tilde{x}(t)= x(-t), \quad \tilde{z}(t) = z(-t), \quad \text{and} \quad \tilde{\theta}(t) = \theta(-t) + \pi.
\end{equation*}
Clearly, $\tilde{\gamma}=(\tilde{x},\tilde{z},\tilde{\theta})$  is a solution of \eqref{ODE system problem} with $\tilde{\gamma}(0)=(x_0,z_0,\pi)$. If $\Sigma$ is neither an entire vertical $\varphi$-graph nor a winglike $\varphi$-minimal surface, then Theorems \ref{teo grafico n simetricos} and  \ref{winglike hiper}, imply  that  on $(0,\infty)$, $\tilde{z}$ is strictly increasing, $\tilde{x}$ is strictly decreasing  and the following inequalities hold
\begin{align*}
& \frac{\pi}{2}<\tilde{\theta} <\pi, \quad 0<\tilde{x}<\infty,\quad \tilde{\theta}'=\dot{\varphi}(\tilde{z}) \cos\tilde{\theta} - \tanh \tilde{x} \sin\tilde{\theta} < 0, \quad \tilde{x}' <0.
\end{align*}
Hence, there exist the limits
$$ \tilde{\theta}_L= \lim_{t\to\infty}\tilde{\theta}(t)\in [ \frac{\pi}{2},\pi),\qquad \tilde{x}_L=\lim_{t\to\infty}\tilde{x}(t)\in [0,x_0),$$
and $  \lim_{t\to\infty}\tilde{x}'(t)=\lim_{t\to\infty}\tilde{\theta}'=0$. Consequently, from system \eqref{ODE system problem}, $\tilde{\theta}_L=\frac{\pi}{2}$, $\tilde{x}_L=0$. Since $\lim_{t\to\infty}\tilde{z}'(t)=1$, it follows that $\lim_{t\to\infty}\tilde{z}(t)=\infty$. Thus, $\Sigma$ must be asymptotic to the plane $\{x=0\}$, which concludes the proof.
\end{proof}
\subsection{Comparing profile curves}
\begin{figure}[H]
    \centering
    \begin{subfigure}[b]{0.45\textwidth}
        \centering
        \includegraphics[width=\textwidth]{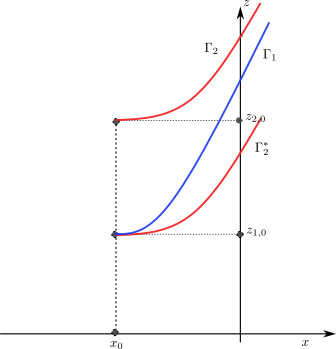}
        \caption{Vertical translation of $\Gamma_2$.}
    \end{subfigure}
    \hfill
    \begin{subfigure}[b]{0.5\textwidth}
        \centering
        \includegraphics[width=\textwidth]{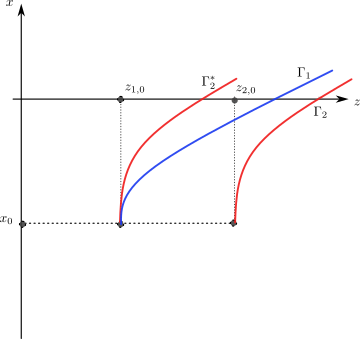}
        \caption{$\Gamma_1$ and $\Gamma_2^*$ as horizontal graph.}
    \end{subfigure}
    \caption{geometric ideas of Proposition \ref{comparar-solucoes}} \label{fig: varia x}
\end{figure}

\begin{figure}[H]
    \centering
    \begin{subfigure}[b]{0.45\textwidth}
        \centering
        \includegraphics[width=\textwidth]{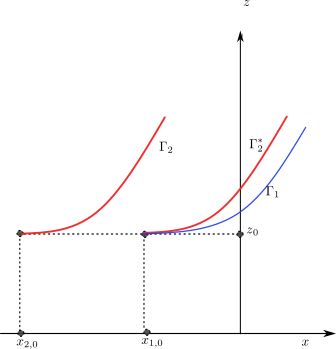}
        \caption{Horizontal translation of $\Gamma_2$.}
    \end{subfigure}
    \hfill
    \begin{subfigure}[b]{0.45\textwidth}
        \centering
        \includegraphics[width=\textwidth]{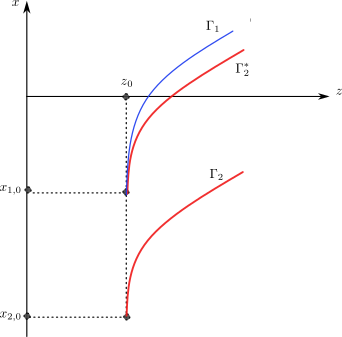}
        \caption{$\Gamma_1$ and $\Gamma_2^*$ as horizontal graph.}
    \end{subfigure}
    \caption{Geometric ideas behind  Proposition \ref{transl vertical curvas}.}\label{fig: varia z}
\end{figure}

We focus next on the geometric study of profile curves whose minimum of the height function occurs either at the same distance to the origin or at the same height.
\begin{prop}\label{comparar-solucoes}
Let $\gamma_1$ and $\gamma_2$ be solutions of \eqref{ODE system problem} satisfying
$$ \gamma_i(0)=(x_{0},z_{i,0},0),
\qquad i=0,1,\qquad
z_{1,0}<z_{2,0}.
$$
Let $ \Gamma_i=\pi\circ\gamma_i $ ($i=1,2$) denote their projections onto the $xz$-plane. Let $\Gamma_2^*$ be the vertical translation of $\Gamma_2$ satisfying that $\Gamma_2^*(0) = \Gamma_1(0)$. If the positive branches  $\Gamma_i^+$ of $\Gamma_i$, $i=1,2$, are  horizontal graphs, then
$\Gamma_1^+$ and $(\Gamma_2^{+})^*$ do not intersect after their common minimum point.
\end{prop}
\begin{proof}
      In a small neighborhood of $x_0$, both curves $\Gamma_1$ and $\Gamma_2^*$ are vertical graphs of smooth functions $u_1$ and $u_2^*$, respectively, such that $u_1$ satisfies \eqref{edo graph} and $u_2^*$ satisfies 
    \begin{equation}
        \left( u^*_2\right)'' = \left(1+     \left(u^*_2\right)^{'^2}\right)\left(\dot{\varphi}\left(u^*_2 + c\right) -\tanh x \left(u^*_2\right)'\right),
    \end{equation}
    with $c = z_{2,0}-z_{1,0}$. 
\\
  Note that $u_2^*(x_0) = u_1(x_0) = z_{1,0}$ and $u_1'(x_0) = u_2^{*'}(x_0) = 0$. Since $z_{2,0}>z_{1,0}$, we have that
    $$ \left(u^*_2\right)''(x_0) = \dot{\varphi}(z_{2,0}) < \dot{\varphi}(z_{1,0}) = u''_1(x_0),$$
and near $x_0$,
    \begin{equation} \label{comparar horizontal x0}
    \left(u^*_2\right) ' < u'_1, \quad \left(u^*_2\right) < u_1.
    \end{equation}
 Now consider these curves as graphs of the horizontal coordinate as a function of height $z$ denoted  by
     $y_1$ and $y_2^*$. The function $y_1$ solves \eqref{graph hori} whereas  $y^*_2$ satisfies 
\begin{equation}\label{graph hori transladado}
        \left(y_2^*\right)'' = \left(1+\left(y^*_2\right)^{'^2}\right)\left(\tanh y^*_2 - \dot{\varphi}\left(z+c \right)\left(y^*_2\right)'\right).
    \end{equation}
    By \eqref{comparar horizontal x0}, in a small neighborhood of $z_{1,0}$,
    $$y^*_2> y_1, \quad \frac{\text{d}}{\text{d}z}y_2^* > \frac{\text{d}}{\text{d}z}y_1.$$
\\
 We claim that $ y^*_2 > y_1$ for all $z > z_{1,0}$. Define the function $d = y_2^* - y_1$. Suppose, by way of contradiction, that there exists $z_m > z_{1,0}$ such that $d'(z_m) = 0$. Then, $z_m$ cannot be a local minimum of $d$. However, evaluating $d''$ at $z_m$ gives
    $$d'' = \left(1 + y_1'\right)\left[\tanh y_2^* - \tanh y_1 - y'_1\left(\dot{\varphi}(z_m + c) -\dot{\varphi}(z_m)\right)\right] > 0,$$
    yielding a contradiction. 
\\  
Therefore, $d'(z)>0$ and consequently $d(z)>0$ for all $z>z_{1,0}$. Thus, $\Gamma_1^+$ and
$\Gamma_2^{*+}$ cannot intersect away from their common minimum
point.
\end{proof}
Figure \ref{fig: varia x} illustrates the geometric idea behind the proof. The left panel shows the vertical translation $\Gamma_2^*$ while the right panel depicts  $\Gamma_2^*$ and $\Gamma_1$ represented as horizontal graphs.  
\begin{theorem}\label{winglikeabovevertical} Let $\gamma_1$ and $\gamma_2$ be solutions of \eqref{ODE system problem} satisfying
$$ \gamma_i(0)=(x_{0},z_{i,0},0),
\qquad i=0,1,\qquad
z_{1,0}<z_{2,0}.
$$
If $\gamma_1$ generates an entire  vertical  $\varphi$-graph, then $\gamma_2$ cannot generate a winglike $\varphi$-surface.
\end{theorem}
\begin{proof}
We may assume that $x_0<0$, otherwise we can compare the left branches. Consider $\Gamma_1$ and $\Gamma_2$ as in  Proposition \ref{comparar-solucoes}. Suppose, for contradiction, that $\Gamma_2$ is winglike. Let $\Gamma_2^*$ be the vertical translation of $\Gamma_2$ satisfying $\Gamma_2^*(0) = \Gamma_1(0)$. Near $x_0$, the curves  $\Gamma_1$ and $\Gamma_2^*$  are vertical graphs of $u_1$ and  $u_2^*$, respectively  with $u_1 > u_2^*$ in a small neighborhood of $x_0$. Because $\Gamma_2^*$ is winglike, $(\Gamma*_2)^+$ it must intersect $\Gamma_1^+$ at some point. However, Proposition \ref{comparar-solucoes} guarantees that $\Gamma_1^+$ and $(\Gamma_2^*)^+$ do not intersect, which is a contradiction.
\end{proof}
Using arguments similar to those in  Theorem \ref{winglikeabovevertical} we obtain the following related results:
\begin{theorem}\label{asymptoticbelowwinglike} Let $\gamma_1$ and $\gamma_2$ be solutions of \eqref{ODE system problem} satisfying
$$ \gamma_i(0)=(x_{0},z_{i,0},0),
\qquad i=0,1,\qquad
z_{1,0}<z_{2,0}.
$$
If $\gamma_1$ generates a vertical $\varphi$-graph of type  (D) in Theorem \ref{withminimum}, then $\gamma_2$ cannot generate a winglike $\varphi$-minimal surface.
\end{theorem}
\begin{theorem}\label{verticalbelowasymptotic} Let $\gamma_1$ and $\gamma_2$ be solutions of \eqref{ODE system problem} satisfying
$$ \gamma_i(0)=(x_{0},z_{i,0},0),
\qquad i=0,1,\qquad
z_{1,0}<z_{2,0}.
$$
If $\gamma_1$ generates an entire vertical  $\varphi$-graph, then $\gamma_2$ cannot generate a vertical $\varphi$-graph of type  (D) in Theorem \ref{withminimum}..
\end{theorem}
We now turn to comparing profile curves whose minimum occurs at the same height.
\begin{prop} \label{transl vertical curvas}
    Let $\gamma_1$ and $\gamma_2$ be solution curve of \eqref{ODE system problem} satisfying
    $$\gamma_i(0) = (x_{i,0}, z_0, 0), \qquad x_{2,0} < x_{1,0}< 0.$$
Let  $
\Gamma_i=\pi\circ\gamma_i$ denote their projections onto the $xz$-plane and let $\Gamma_2^*$ be the horizontal translation of $\Gamma_2$ such that $\Gamma_2^*(0) = \Gamma_1(0)$. If the positive branches  $\Gamma_i^+$ of $\Gamma_i$, $i=1,2$ are  horizontal graphs, then $(\Gamma_2^+)^*$ does not intersect $\Gamma_1^+$.
\end{prop}
\begin{proof}
    First, observe that near $x_0$, $\Gamma_2^*$ is a vertical graph of a function $u_2^*$ satisfying
    \begin{equation} \label{graph vertical transladado}
      (u_2^*)'' = \bigl(1 + (u_2^*)^{'2}\bigr)\Bigl(\dot{\varphi}(u_2^*) - (u_2^*)' \tanh(x - c_0)\Bigr),  
    \end{equation}
where $c_0 = x_{1,0} - x_{2,0}$, while $\Gamma_1$ is a vertical graph of $u_1$ which solve  \eqref{edo graph}. For $\varepsilon > 0$, define $d: [x_0, x_0 + \varepsilon) \to \mathbb{R}$ by
$$ d(x) = u_2^*(x) - u_1(x). $$
Note that $d(x_0) = 0$ and $d'(x_0) = 0$, since $u^*_2(x_0) = u_1(x_0) = z_0$ and $u^{'*}_2(x_0) = u'_1(x_0) = 0$. Moreover, 
$$ d''(x_0) = \dot{\varphi}\left(u^*_2(x_0)\right) - \dot{\varphi}\left(u_1(x_0)\right) = 0.$$
However, direct computation yields
$$d'''(x_0) =  - \dot{\varphi}(z_0) \left( \tanh(x_0 - c_0) + \tanh(x_0) \right) > 0.$$
Hence, in a small neighborhood $(x_0, x_0 + \varepsilon)$ of $x_0$, we have 
\begin{equation}\label{d>0ed'>0}
    d > 0, \quad d' > 0. 
\end{equation}
Now consider these curves as horizontal graphs. In this setting, $\Gamma_1$ is the horizontal graph of a function $y_1$ satisfying \eqref{graph hori}, and $\Gamma_2^*$ is a horizontal graph of a function $y_2^*$ satisfying
$$ (y_2^*)'' = \bigl(1 + (y_2^*)^{'2}\bigr)\Bigl(\tanh(y_2^* - c_0) - (y_2^*)' \dot{\varphi}(z)\Bigr). $$
Furthermore, by \eqref{d>0ed'>0},  $y_2^* < y_1$ and $(y_2^*)' < y_1'$, near $z_0$ . Define
$ D(z) = y_1(z) - y_2^*(z). $
We will show that $D(z) > 0$ for all $z > z_0$. Since $D > 0$ and $D' > 0$ in a small neighborhood of $z_0$, it is sufficient to show that $D'(z) > 0$ for all $z > z_0$. Suppose, for contradiction, that  $D'(z_m) = 0$ for some $z_m>z_0$. Taking $z_m$ to be the first critical point of $D$, we must have $D''(z_m) \leq 0$. However,
\begin{equation} \label{seg dev d}
   D''(z_m) = \bigl(1 + (y_1'(z_m))^2\bigr)\Bigl(\tanh(y_1\left(z_m\right)) - \tanh\left(y^*_2(z_m) - c_0\right)\Bigr) > 0,  
\end{equation}
which is a contradiction.
\end{proof}
Figure \ref{fig: varia z} illustrates the geometric idea behind the proof. In the next proposition, we prove that no winglike $\varphi$ surfaces can lie to the left of a vertical graph.
\begin{theorem} \label{winglikeleftgraph} If 
   Let $\gamma_1$ and $\gamma_2$ be solutions of \eqref{ODE system problem} satisfying
$$ \gamma_i(0)=(x_{i,0},z_{0},0),
\qquad i=0,1,\qquad
x_{2,0}<x_{1,0}<0.
$$
If $\gamma_1$ generates a  winglike $\varphi$-minimal surface, then $\gamma_2$ cannot generate an entire  vertical  $\varphi$-graph.
\end{theorem}
\begin{proof} 
Suppose, for contradiction, that the solution curve $\gamma_2$ of \eqref{ODE system problem} with initial condition $\gamma_2(0) = (x_{2,0}, z_0, 0)$  satisfies that $\Gamma_2$ is a vertical graph. Consider $\Gamma_2^*$, the horizontal translation of $\Gamma_2$, satisfying $\Gamma_2^*(0) = \Gamma_1(0)$. 
\\
Since $\Gamma_2^*$ is a vertical graph and $\Gamma_1$ is a winglike curve,  $(\Gamma_2^*)^+$  must intersect $\Gamma_1^+$. However,  Proposition \ref{transl vertical curvas} asserts that $(\Gamma_2^*)^+$ does not intersect $(\Gamma_1)^+$. This contradiction completes the proof.
\end{proof}
By reasoning as in Theorem \ref{winglikeleftgraph} we also obtain, 
\begin{theorem} \label{winglikeleftasymptotic} 
   Let $\gamma_1$ and $\gamma_2$ be solutions of \eqref{ODE system problem} satisfying
$$ \gamma_i(0)=(x_{i,0},z_{0},0),
\qquad i=0,1,\qquad
x_{2,0}<x_{1,0}<0.
$$
If $\gamma_1$ generates a winglike $\varphi$-minimal surface, then $\gamma_2$ cannot generate a vertical  $\varphi$-graph lof type (D) in Theorem \ref{withminimum}.
\end{theorem}
\begin{theorem} \label{asymptoticleftgraph} If 
   Let $\gamma_1$ and $\gamma_2$ be solutions of \eqref{ODE system problem} satisfying
$$ \gamma_i(0)=(x_{i,0},z_{0},0),
\qquad i=0,1,\qquad
x_{2,0}<x_{1,0}<0.
$$
If $\gamma_1$ generates a vertical  $\varphi$-graph of type (D) in Theorem \ref{withminimum}, then $\gamma_2$ cannot generate an entire vertical  $\varphi$-graph.
\end{theorem}
\subsection{Limits of profile curves}
\begin{figure}[ht]
    \centering
    \includegraphics[scale=0.7]{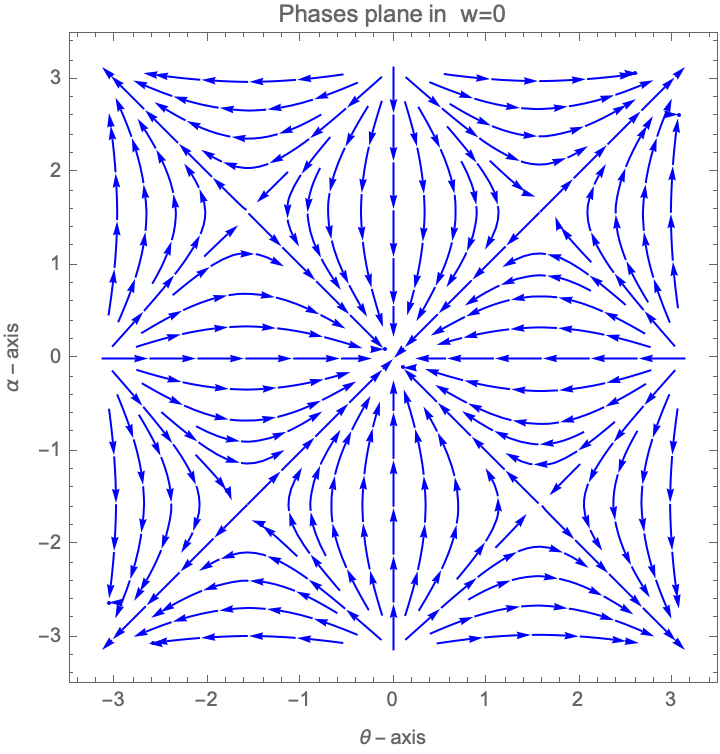}
    \caption{Trajectories of \eqref{system w=0} }
    \label{fasesw0}
\end{figure}

Now, we analyze which solutions of    \eqref{ODE system problem} converge to the solution without a minimum obtained in Theorems \ref{existence} and  \ref{sem minimo teorema}.
\begin{theorem}\label{limittoa} Let $x_0\neq0$. Then there exists $z_0>a$ such that for any $z_1\in (a,z_0)$,  the solution $\gamma_{z_1}$ of \eqref{ODE system problem} with $\gamma_{z_1}(0)=(x_0,z_1,0)$ is of winglike type.
\end{theorem}
\begin{proof}
Let $x=y_{x_0}(z)$, $z>a$ be the solution of \eqref{equiv ode zero} given by Theorems \ref{existence} and  \ref{sem minimo teorema} and let $\gamma_{z_0}=(x,z,\theta)$ be the solution of  \eqref{ODE system problem} satisfying $\gamma_{z_0}(0)=(x_0,z_0,0)$, $x_0>0$. Let $\Gamma^+_{z_0}=(x^+,z^+)=(x,z)|_{(0,\infty)}$ and $\Gamma^-_{z_0}=(x^-,z^-)=(x,z)|_{(-\infty,0)}$ denote  the right and left branches of $(x,z)$, respectively. By Theorems \ref{existence} and \ref{sem minimo teorema}, $\Gamma^+_{z_0}$ and $\Gamma^-_{z_0}$ are horizontal graphs given by  $x^+ = y^+_{z_0}(z)$ and  $x^- = y^-_{z_0}(z)$, $z> z_0$ for  some functions $y^\pm_{z_0}$ that solve \eqref{graph hori} on $(z_0,\infty)$ and  $y^{\pm}_{z_0}(z_0)=x_0$.\\
If $\lim_{z_0\to a} y^{\pm}_{z_0} = y^{\pm}_{x_0}$, then, arguing as in Lemma \ref{initial derivative zero}, we obtain that $y^{\pm}_{x_0}$ are solutions of \eqref{equiv ode zero}.
 By the uniqueness of the solutions to the problem \eqref{equiv ode zero}, we obtain 
\begin{equation}\label{branches}
\lim_{z_0\to a} y^+_{x_0}= y_{x_0}, \qquad \lim_{z_0\to a} y^-_{x_0}= y_{x_0}.
\end{equation}
However, $y_{x_0}(z)>x_0$, for all $z>a$. Thus, if  $\gamma_{z_0}$ were of type (A), (B) or (D) in Theorem \ref{withminimum} as $z_0\to a$, limit \eqref{branches} would lead to a contradiction for the left branch $\gamma^-_{z_0}$ of $\gamma_{z_0}$. Consequently, there exists $z_0$ close to $a$ such that $\gamma_{z_0}$ is winglike. The desired conclusion then follows by applying   Theorems \ref{winglikeabovevertical} and \ref{asymptoticbelowwinglike}.
\end{proof}
To study the trajectories of  \eqref{ODE system problem}  near  the boundary of its phase space, we perform the change of variables, 
$$ w =\frac{1}{z-a}, \qquad \cos \alpha = \tanh x.$$
In terms of $(\alpha,w,\theta)$  the system  \eqref{ODE system problem}   becomes
\begin{equation}\label{new system}\left\{
    \begin{array}{ll}
     w' =& - w^2 \sin \theta,\\
         \alpha'=& - \cos \theta \sin \alpha  \\
         \theta' =&  \psi(w)\cos\theta -\sin \theta \cos \alpha 
    \end{array}\right.
\end{equation}
where $\psi(w)$ is defined by
\begin{equation*}\psi(w)=\left\{
    \begin{array}{ll}
   \dot{\varphi}(1/w), &w>0,\\
    g(w), \quad &w\leq 0
    \end{array}\right.
\end{equation*}
for a fixed function $g\in{\cal C}^1(\mathbb{R})$  satisfying
$$ g(0)=0, \qquad g'(0)=L_\infty= \lim_{z\to\infty}z^2 \ddot{\varphi}(z).$$
Under hypothesis (H2), the system \eqref{new system} is well defined and admits unique solutions for any initial conditions in $\mathbb{R}^3$. Moreover any solution $\gamma=(x,z,\theta)$  of \eqref{ODE system problem} yields a solution $\hat{\gamma}=(w,\alpha,\theta)$ of \eqref{new system}, via $w (z-a)=1$ and $\cos\alpha = \tanh x$.
\begin{theorem}\label{largeheight} If $x_0\neq 0$, then there exists $z_0>a$, large enough such that for any  $z_1>z_0$,  the solution $\gamma_{z_1}$ of \eqref{ODE system problem} satisfying $\gamma_{z_1}(0)=(x_0,z_1,0)$  generates an entire vertical $\varphi$-graph.
\end{theorem}
\begin{proof}
For any $z_0>a$, let  $\gamma_{z_0}=(x,z,\theta)$ be the solution of  \eqref{ODE system problem} with
 $\gamma_{z_0}(0)=(x_0,z_0,0)$. Then $\hat{\gamma}_{w_0}=(w,\alpha,\theta)$ is the solution of \eqref{new system} satisfying 
 $$\hat{\gamma}_{w_0}(0)=(w_0,\alpha_0,0), \quad w_0=\frac{1}{z_0-a}, \quad \cos \alpha_0=\tanh x_0, \quad \alpha_0\in (0,\pi).$$
 Moreover, $\lim_{w_0\to 0} \hat{\gamma}_{w_0}=\hat{\gamma}^*$, where $\hat{\gamma}^*$ is the solution of \eqref{new system} passing through $(0,\alpha_0,0)$.
 
 Since $\lim_{w\to 0} \psi(w) =0$,  the plane $\Pi\equiv \{w=0\}$ is invariant  under  system \eqref{new system}. For any initial condition $(0,\alpha_0,\theta_0)$ there exists a unique solution $\hat{\gamma}=(0, \alpha,\theta)$ of 
 \begin{equation}\label{system w=0}\left\{
    \begin{array}{ll}
         \alpha'=& - \cos \theta \sin \alpha  \\
         \theta' =&   -\sin \theta \cos \alpha,
    \end{array}\right.
\end{equation}
passing through $(0,\alpha_0,\theta_0)$.  In particular, for the initial condition $(0,\alpha_0,0)$, $\alpha_0<\frac{\pi}{2}$,  the solution of \eqref{new system} through this point is $\hat{\gamma}^*(t)=(0, \alpha(t),0)$, with $\cos \alpha(t) = \tan\frac{\alpha_0}{2} \ {\rm e}^{-t} $. As shown in Figure \ref{fasesw0},  this solution satisfies 
$$ \lim_{t\to -\infty} \alpha(t)= \pi, \qquad  \lim_{t\to +\infty} \alpha(t)= 0,$$
and $\cos \alpha$ changes sign when $t$ is sufficiently close to $\log \tan\frac{\alpha_0}{2}$. The case $\alpha_0>\frac{\pi}{2}$  is treated analogously.

Therefore, for $z_0$ large enough, the function $x(t)$ of $\gamma_{z_0}$ must vanish at some point and Theorem \ref{withminimum} then implies that,
$$ \lim_{t\to -\infty} x(t)= -\infty (+\infty), \qquad  \lim_{t\to +\infty} x(t)= +\infty(-\infty).$$
That is, $\gamma_{z_0}$ generates an entire vertical $\varphi$-graph. The proof is completed by applying  Theorems \ref{winglikeabovevertical} and \ref{verticalbelowasymptotic}.
\end{proof}
\begin{remark} If $\theta_0\neq 0$,   Figure \ref{fasesw0} shows that  for any sufficiently large $z_0$, there exist $z_1>z_0$ and $x_0\neq 0$ such that the solution $\gamma $ of  \eqref{ODE system problem} with $\gamma(0) = (x_0,z_1,\theta_0)$ is winglike.
\end{remark}
Now, we examine solutions $\gamma$ of  \eqref{ODE system problem} with $\gamma(0) = (x_0,z_0,0)$ when  $|x_0|\to\infty$  (i.e., as $\alpha_0\to 0$ or $\alpha_0\to \pi$). First, note that the planes $\Pi_{-}1\equiv\{\alpha=\pi\}$ and $\Pi_1\equiv\{\alpha=0\}$ are invariants under the system \eqref{new system}. Considering $\Pi_1$ (with $\Pi_{-1}$ handled similarly) if  $\gamma=(\alpha,w,\theta)$ is a solution of \eqref{new system} with $\gamma(0)=(0,w_0,\theta_0)$ then $\alpha\equiv 0$ and the system  \eqref{new system}  reduces to
\begin{equation}\label{system alpha=0}\left\{
    \begin{array}{ll}
         w' =& - w^2 \sin \theta  \\
         \theta' =& \psi(w) \cos \theta  -\sin \theta,
    \end{array}\right.
\end{equation}
\begin{figure}[ht]
    \centering
    \includegraphics[scale=0.7]{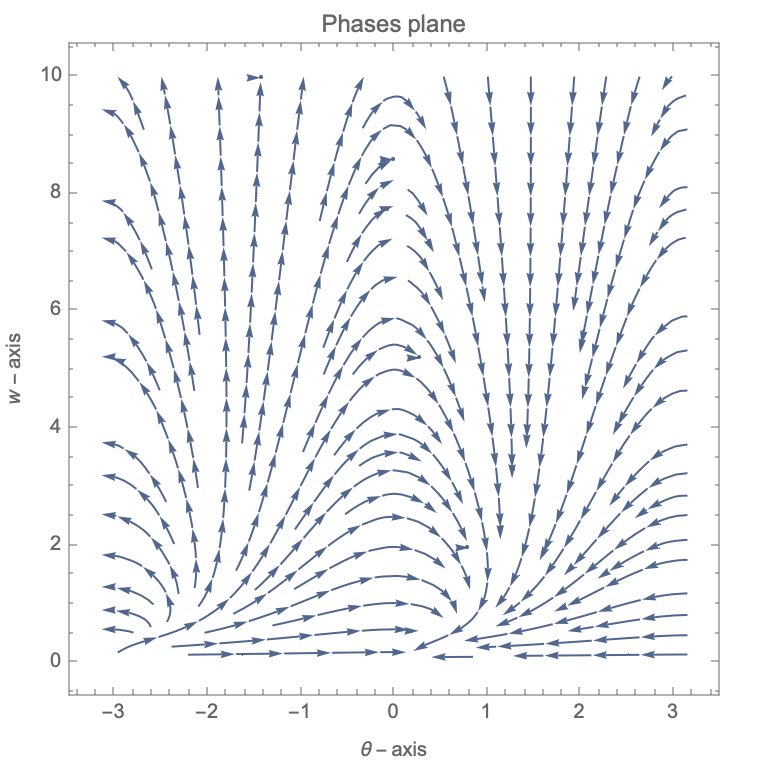}
    \caption{Trajectories of \eqref{system alpha=0} }
    \label{fasesalpha0}
\end{figure}
\begin{lemma}\label{sol alpha0} Let $(w,\theta)$ be the solution of \eqref{system alpha=0} with $w(0)=w_0>0$ and $\theta(0)=0$, then $w$  and $\theta$ are defined for all $t\in\mathbb{R}$, both possess a unique critical point (which is a local  maximum) and (see Figure \ref{fasesalpha0}),
\begin{align*}
\lim_{t\to -\infty} w(t) =0, & \lim_{t\to +\infty} w(t) =0,\\
\lim_{t\to -\infty} \theta(t) =-\pi, & \lim_{t\to +\infty} \theta(t) =0.
\end{align*}
\end{lemma}
\begin{proof} Since the plane $\Pi\equiv\{w=0\}$ is invariant under \eqref{new system}, we have  $w>0$ on its maximal interval of definition $I_M=(t_{min},t_{max})$. If $w'(t_0)=0$, for some $t_0\in I_M$, then $$w''(t_0)= - \psi(w(t_0))\cos^2\theta(t_0) w(t_0)<0,$$ so $w$ has a unique critical point at $t=0$ which is a maximum.
Consequently, $0<w<w_0$ on $I_M$. Following  Proposition 4.2 of \cite{freitas2026geometry} it follows that $I_M=\mathbb{R}$, meaning the solution is globally defined.
\\
If $\theta'(t_0)=0$, then $$\theta''(t_0)=\frac{\ddot{\varphi}(1/w)}{w\dot{\varphi}(1/w)}\sin^2\theta<0,$$ so $\theta$ has at most a critical point which is also a maximum. 
\\
Due to  the boundedness and monotonicity of $w$ and $\theta$ near infinity, system \eqref{system alpha=0} yields 
$$0 = -\lim_{t\to \pm\infty}\sin \theta(t), \quad  0=\lim_{t\to \pm\infty}w(t)$$
and  $\theta'$ changes  sign on $(0,\infty)$ completing the proof.
\end{proof}
\begin{theorem} \label{alphato0}If  $z_0>a$, then there exists $x_w<0$ (respectively, $x_w>0$) such that for any  $x_1<x_w$ (respectively, $x_1>x_w$) the solution $\gamma_{x_1}$ of  \eqref{ODE system problem} with $\gamma_{x_1}(0) = (x_1,z_0,0)$ generates a winglike $\varphi$-minimal surface.
\end{theorem}
\begin{proof} Assume $x_0<0$ without loss of generality. Let $\gamma_{x_0}=(x,z,\theta)$ solve \eqref{ODE system problem} with $\gamma_{x_0}(0)=(x_0,z_0,0)$. Then $\hat{\gamma}_{\alpha_0}=(\alpha,w,\theta)$ solves \eqref{new system} with initial value $\hat{\gamma}_{\alpha_0}=(\alpha_0,w_0,0)$, where $\cos \alpha_0=\tanh x_0$, $w_0 (z_0-a)=1$. Since $\Pi_1$ is invariant under \eqref{new system},
$$\lim_{\alpha_0\to 0}\hat{\gamma}_{\alpha_0} = \tilde{\gamma}_{w_0} = (0, \tilde{w},\tilde{\theta})$$
where $( \tilde{w},\tilde{\theta})$ is the solution of \eqref{system alpha=0} with $ \tilde{w}(0)=w_0$ and $\tilde{\theta}(0)=0$. 

Thus, $\alpha$ remains close to $0$ whenever $\alpha_0$ is chosen sufficiently small, that is, there exists $x_0$ with $|x_0|$ sufficiently  large such that  $x(t)<M<0$ on $(-\infty,\infty)$. By Theorem \ref{withminimum},  $\gamma_{x_0}$ must generate a winglike $\varphi$-minimal surface. The conclusion follows from Theorems \ref{winglikeleftgraph} and  \ref{winglikeleftasymptotic}.
\end{proof}
\begin{theorem}\label{alphatopi2}
If  $z_0>a$, then there exists $x_w<0$ (respectively, $x_g>0$) such that for any  $x_1\in(x_g,0)$ (respectively, $x_1\in (0,x_g)$) the solution $\gamma_{x_1}$ of  \eqref{ODE system problem} with $\gamma_{x_1}(0) = (x_1,z_0,0)$ generates a non-symmetric entire vertical  $\varphi$-graph.
\end{theorem}
\begin{proof} Assume $x_0<0$ and let $\gamma_{x_0}=(x,z,\theta)$ be the solution of \eqref{ODE system problem} with $\gamma_{x_0}(0)=(x_0,z_0,0)$. By Theorem \ref{bowl hiperbolico} and uniqueness of  solutions, $\lim_{x_0\to 0}\gamma_{x_0}$ is the symmetric entire vertical graph with minimum at $(0,z_0)$, generated by  $\tilde{\gamma}_{z_0}=(\tilde{x},\tilde{z},\tilde{\theta})$ with $\tilde{\gamma}_{z_0}(0)=(0,z_0,0)$. Since $\tilde{x}$ changes sign at $t=0$, it follows that, for $x_0$ sufficiently close to $0$ , $x$ also changes  sign. Thus, $\gamma_{x_0}$ generates a non-symmetric entire vertical graph and the assertion follows from Theorems \ref{winglikeleftgraph},  \ref{winglikeleftasymptotic} and \ref{verticalbelowasymptotic}.
\end{proof}

\subsection{Existence of vertical graphs in a half-space of \texorpdfstring{$\mathbb{H}^2\times \mathbb{R}$}{HR}}
We now turn to the existence of the most significant case of hyperbolic $\varphi$-minimal surface: vertical graphs restricted to a  half-space  $\{x<0\}$ (or  $\{x>0\}$)  of $\mathbb{H}^2\times \mathbb{R}$.  Unlike the symmetric case and the winglike type surfaces, which are well-documented in the literature on product spaces (see  \cite{de2019translating} and \cite{lima-pipoli}) and exhibit behavior analogous to  rotationally symmetric and parabolic examples,  no other vertical graph of this type exists in the rotational and parabolic cases. A surface with comparable geometric features appears in \cite{lima-pipoli} and \cite{bueno2024new}  as a special case of a translating soliton.

\

\begin{figure}[htp!]
    \centering
    \begin{subfigure}[b]{0.45\textwidth}
        \centering
        \includegraphics[width=\textwidth]{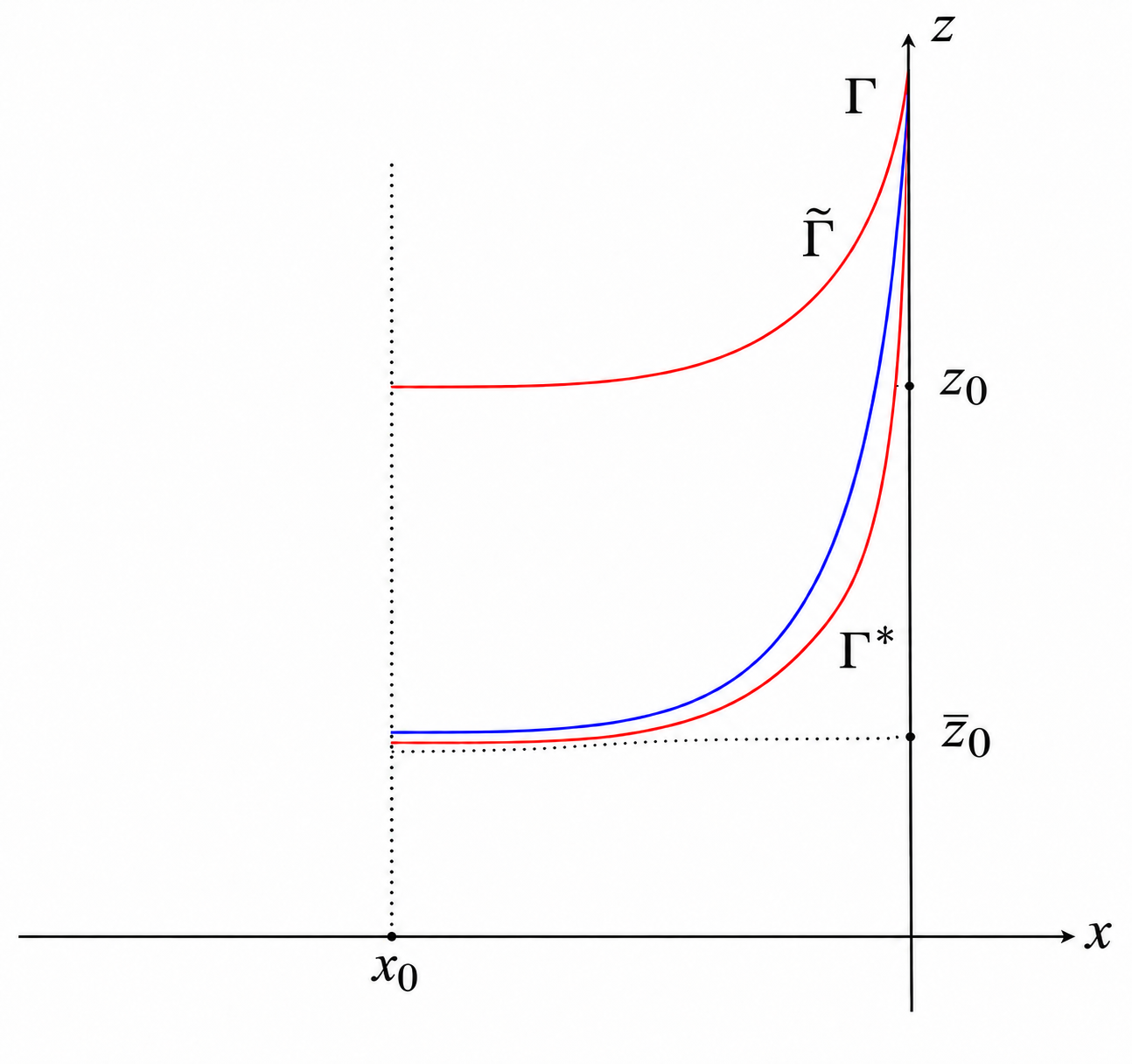}
        \caption{Vertical translation of $\bar{\Gamma}$.}
    \end{subfigure}
    \hfill
    \begin{subfigure}[b]{0.45\textwidth}
        \centering
        \includegraphics[width=\textwidth]{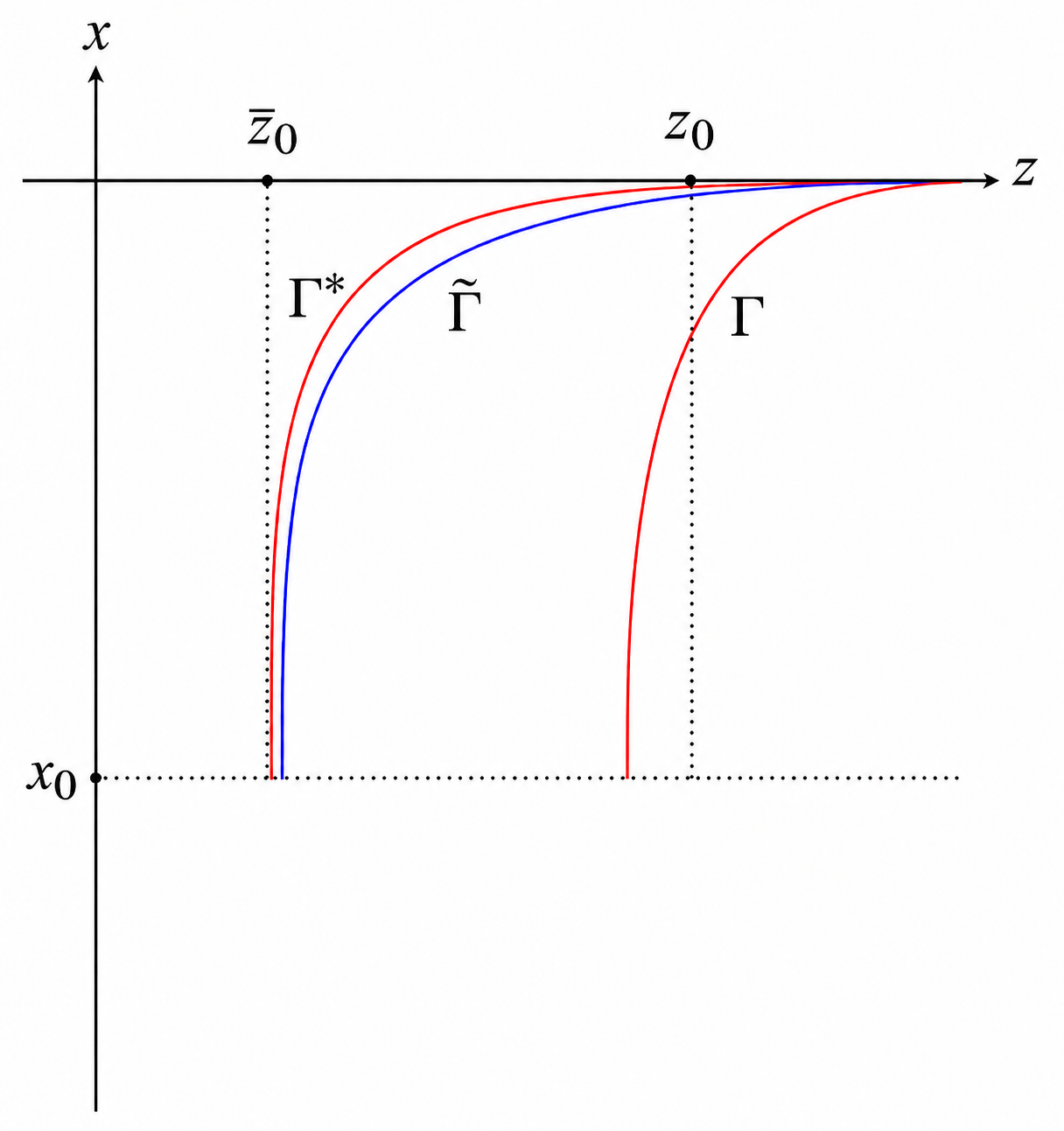}
        \caption{$\Gamma$ and $\bar{\Gamma}*$ as horizontal graph.}
    \end{subfigure}
    \caption{Geometric ideas of the proof of Proposition \ref{existenceasymptotic}.}\label{varia x} 
\end{figure}

For a fixed  $x_0\neq 0$  and  each $z_0>a$, let $\gamma_{z_0}$ denote the solution of \eqref{ODE system problem} with $\gamma_{z_0}(0)=(x_0,z_0,0)$. We define,
\begin{align*}
W_{x_0} & = \{ z>a \ | \ \gamma_z \text{ generates a winglike $\varphi$-minimal surface}\},\\
G_{x_0} & = \{ z>a \ | \ \gamma_z \text{ generates a non-symmetric entire vertical  $\varphi$-graph}\}.
\end{align*}
By Theorems \ref{winglikeabovevertical}, \ref{limittoa} and \ref{largeheight}, $W_{x_0}$ is bounded from above, $G_{x_0}$ is bounded from below and $$\sup W_{x_0} \leq \inf G_{x_0}.$$
More precisely, we establish the following threshold result:
\begin{prop}\label{existenceasymptotic} There exists a constant $b>a$, such that $W_{x_0}=(a,b)$ and $G_{x_0}=(b,\infty)$. Moreover $\gamma_b$ generates a vertical $\varphi$-graph in the half-space $\{x>0\}$ (or $\{x<0\})$ asymptotic to the vertical plane $\{x=0\}$ with an unbounded height function  at infinity.
\end{prop}
\begin{figure}[ht!]
    \centering
\includegraphics[width=0.5\linewidth]{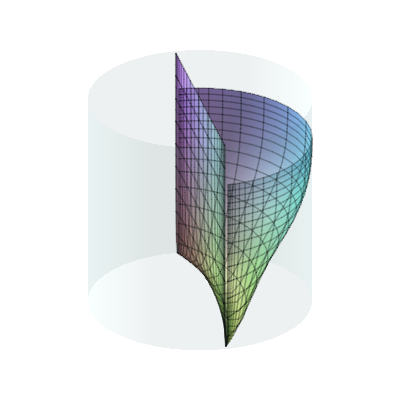}
    \caption{A graph over a half-space for $\dot{\varphi}(z) = \frac{1}{z}$}
\end{figure}
\begin{proof}
First, we show that $W_{x_0}$ is open. Let $z_0\in W_{x_0}$.  Since  $\gamma_{z_0} =(x,z,\theta)$ is winglike,  Theorem \ref{winglike hiper} implies that there exists $t_1\in \mathbb{R}$ such that $\theta(t_1)=\pm\frac{\pi}{2}$. Continuous dependence on initial conditions guarantees, that  for any $\bar{z}_0$ sufficiently close to $z_0$, the solution $\gamma_{{\bar z}_0}=(\bar x,\bar z,\bar \theta)$ satisfies $\bar\theta(\bar{t}_1)=\pm\frac{\pi}{2}$ for $\bar{t}_1$ near $t_1$.
\\
Similarly, if $z_0\in G_{x_0}$,  Theorem \ref{teo grafico n simetricos} implies that  if  $\gamma_{z_0} =(x,z,\theta)$  there exists $t_1\in \mathbb{R}$ such that $x(t_1)=0$ and continuous dependence on  initial conditions again shows that $G_{x_0}$ is open. 
\\
By Theorems \ref{winglikeabovevertical}, \ref{asymptoticbelowwinglike} and \ref{verticalbelowasymptotic}, there exist $b,c\in \mathbb{R}$ with  $a<b\leq c$ such that $W_{x_0}=(a,b)$ and $G_{x_0}=(c,\infty)$. Assume for contradiction that $b<c$. Then   the solutions $\gamma_c=(x,z,\theta)$ and $\gamma_b=(\bar x,\bar z,\bar\theta)$ generate $\varphi$-minimal  vertical graph in the half-space $\{x<0\}$ asymptotic to $x=0$. Let $\Gamma$ and  $\bar{\Gamma}$ be the projections of $\gamma_c$ and $\gamma_b$ onto the $xz$-plane.  As in Proposition \ref{comparar-solucoes}, let $\Gamma^*$ be the vertical translation of $\Gamma$ with $\Gamma^*(0) = (x_0, b)$. 
If $(\Gamma^*)^+$ intersects $(\bar{\Gamma})^+$,  Proposition \ref{comparar-solucoes} yields a contradiction, showing that no vertical graph in the half-space $\{x<0\}$ (or $\{x>0\}$) can lie above another.
   \\
    Now suppose $(\Gamma^*)^+$ does not intersect $(\bar{\Gamma})^+$. Represent $(\bar{\Gamma})^+$ as a horizontal graph of $\bar{y}$ and $(\bar{\Gamma}^*)^+$ as a horizontal graph of $y^*$. Here, $\bar{y}$  solves \eqref{graph hori} while  $y^*$ solves \eqref{graph hori transladado}. As established in Proposition \ref{comparar-solucoes}, 
    $d(z) = y^*(z) - \bar{y}(z)$
    satisfies     $ d > 0$  and $ d' > 0$ in a small neighborhood of $\bar{z}_0$.
    Furthermore, by the asymptotic properties of  half-space graphs  in  $\mathbb{H}^2\times \mathbb{R}$ we have that 
    $$\lim_{z \to \infty} d(z) = 0.$$
    Hence, there must exists $z_m > z_0$ such that $d'(z_m) = 0$. Taking  $z_m$ to be the first critical point of $d$ we must have $d''(z_m) \leq 0$. However,  equation \eqref{seg dev d} shows that  $d'' > 0$, yielding a contradiction. Thus, $b=c$ completing the proof (see Figure \ref{varia x} for some geometric ideas of this proof).
\end{proof}

By applying an analogous argument for fixed  $z_0>a$ and varying $x_0<0$, let $\gamma_{x_0}$ be the solution of \eqref{ODE system problem} with $\gamma_{x_0}(0)=(x_0,z_0,0)$. Define
\begin{align*}
W_{z_0} & = \{ x<0 \ | \ \gamma_x \text{ generates a winglike $\varphi$-minimal surface}\},\\
G_{z_0} & = \{ x<0 \ | \ \gamma_x \text{ generates a non-symmetric entire vertical  $\varphi$-graph}\}.
\end{align*}
By Theorems \ref{winglikeabovevertical}, \ref{alphato0} and \ref{alphatopi2}, $W_{z_0}$ is bounded from above, $G_{z_0}$ is bounded from below and $\sup W_{x_0} \leq \inf G_{x_0}.$ This leads to:
\begin{prop}\label{existenceasyntotic2} There exists $b$, $b<0$, such that $W_{z_0}=(-\infty,b)$ and $G_{z_0}=(b,0)$. Moreover $\gamma_b$ generates a vertical $\varphi$-graph in the half-space $\{x<0\}$, asymptotic to the vertical plane $\{x=0\}$, and with unbounded height function  at infinity.
\end{prop}
\subsection{Classification}
From Theorems \ref{unicas} and \ref{withminimum} we can obtain the  following classification of all $\varphi$-minimal surfaces invariant under hyperbolic translation,
\begin{theorem}
    Let $\Sigma$ be a hyperbolic $\varphi$-minimal surface in $\mathbb{H}^2\times \mathbb{R}$. If $\Sigma$ achieves a point of minimal height, then $\Sigma$ lies in either a entire vertical graph (symmetric or not), a vertical graph in the half-space $\{x<0\}$ (or $\{x>0\}$), or a winglike surface. 
    If $\Sigma$ does not have a minimum point, then $\Sigma$ lies  in one of the surfaces described in Theorem \ref{sem minimo teorema}.
\end{theorem}

\subsubsection*{Acknowledgments and funding.} Katianny Rosa Freitas was supported by CAPES Finance Code 001 and she thanks the Departamento de Geometría y Topología at Universidad de Granada for the hospitality, where part of this work was conducted. João Paulo dos Santos was partially supported by CNPq 313200/2025-4 and FAPDF 0019300001678/2024-39.  Antonio Martínez was partially supported by Ministerio de Ciencia,  Innovaci\'on y Universidades
Grant No: PID2024-160586NB-100; “Maria de Maeztu” Excellence Unit IMAG, reference CEX2020-
001105-M and Junta de Andaluc\'ia Grant No. FQM325.

\newpage
\bibliographystyle{plain}
\bibliography{sample}

@article{lima-pipoli,
    AUTHOR = {de Lima, Ronaldo F. and Pipoli, Giuseppe},
     TITLE = {Translators to higher order mean curvature flows in {$\mathbb{R}^n \times\mathbb{ R}$} and {$\mathbb{ H}^n\times\mathbb{R}$}},
   JOURNAL = {J. Geom. Anal.},
  FJOURNAL = {Journal of Geometric Analysis},
    VOLUME = {35},
      YEAR = {2025},
    NUMBER = {3},
     PAGES = {Paper No. 92, 50},
      ISSN = {1050-6926,1559-002X},
   MRCLASS = {53E40 (53C42 53C45)},
  MRNUMBER = {4863897},
MRREVIEWER = {Peter\ Olamide\ Olanipekun},
       DOI = {10.1007/s12220-025-01925-5},
       URL = {https://doi.org/10.1007/s12220-025-01925-5},
}

@article{bueno2024new,
    AUTHOR = {Bueno, Antonio and L\'{o}pez, Rafael},
     TITLE = {A new family of translating solitons in hyperbolic space},
   JOURNAL = {J. Korean Math. Soc.},
  FJOURNAL = {Journal of the Korean Mathematical Society},
    VOLUME = {63},
      YEAR = {2026},
    NUMBER = {2},
     PAGES = {295--311},
      ISSN = {0304-9914,2234-3008},
   MRCLASS = {53A10 (53C21 53C42 53E10)},
  MRNUMBER = {5037056},
       DOI = {10.4134/JKMS.j250079},
       URL = {https://doi.org/10.4134/JKMS.j250079},
}

@article{martinez2022equilibrium,
  title={Equilibrium of surfaces in a vertical force field},
  author={Mart{\'\i}nez, Antonio and Mart{\'\i}nez-Trivi{\~n}o, AL},
  journal={Mediterranean Journal of Mathematics},
  volume={19},
  number={1},
  pages={1--28},
  year={2022},
  publisher={Springer}
}

@article{bueno2024classgrimreapersmathbbh2timesmathbbr,
    AUTHOR = {Bueno, Antonio and L\'{o}pez, Rafael},
     TITLE = {The class of grim reapers in {$\mathbb{ H}^2\times \mathbb{R}$}},
   JOURNAL = {J. Math. Anal. Appl.},
  FJOURNAL = {Journal of Mathematical Analysis and Applications},
    VOLUME = {541},
      YEAR = {2025},
    NUMBER = {2},
     PAGES = {Paper No. 128730, 23},
      ISSN = {0022-247X,1096-0813},
   MRCLASS = {53A05 (53A35 53E10)},
  MRNUMBER = {4783578},
MRREVIEWER = {Zhi\ Li},
       DOI = {10.1016/j.jmaa.2024.128730},
       URL = {https://doi.org/10.1016/j.jmaa.2024.128730},
}

@article{de2019translating,
  title={Translating solitons in Riemannian products},
  author={de Lira, Jorge HS and Martín, Francisco},
  journal={Journal of Differential Equations},
  volume={266},
  number={12},
  pages={7780--7812},
  year={2019},
  publisher={Elsevier}
}

@article{bueno2018translating,
  title={Translating solitons of the mean curvature flow in the space $\mathbb{H}^2\times\mathbb{R}$},
  author={Bueno, Antonio},
  journal={Journal of Geometry},
  volume={109},
  pages={1--26},
  year={2018},
  publisher={Springer}
}

@article{de2020hypersurfaces,
    AUTHOR = {de Lima, R. F. and Manfio, F. and dos Santos, J. P.},
     TITLE = {Hypersurfaces of constant higher-order mean curvature in
              {$M\times\mathbb{R}$}},
   JOURNAL = {Ann. Mat. Pura Appl. (4)},
  FJOURNAL = {Annali di Matematica Pura ed Applicata. Series IV},
    VOLUME = {201},
      YEAR = {2022},
    NUMBER = {6},
     PAGES = {2979--3028},
      ISSN = {0373-3114,1618-1891},
   MRCLASS = {53B25 (53C24 53C42)},
  MRNUMBER = {4498418},
       DOI = {10.1007/s10231-022-01229-3},
       URL = {https://doi.org/10.1007/s10231-022-01229-3},
}

@article{freitas2026geometry,
  title={Geometry of a class of {$\varphi$}-minimal surfaces in {$\mathbb{H}^2 \times \mathbb{R}$}},
  author={Freitas, Katianny Rosa and Santos, Jo{\~a}o Paulo dos and Mart{\'\i}nez, Ant{\'o}nio},
  journal={arXiv preprint arXiv:2610.27187},
  year={2026}
}

@article{mean-convex,
    AUTHOR = {Mart\'{\i}nez, Antonio and Mart\'{\i}nez-Trivi\~{n}o, Antonio
              Luis and dos Santos, Jo\~{a}o Paulo},
     TITLE = {Mean convex properly embedded
              {$[\varphi, \vec{e}_3]$}-minimal surfaces in
              {$\mathbb{R}^3$}},
   JOURNAL = {Rev. Mat. Iberoam.},
  FJOURNAL = {Revista Matem\'{a}tica Iberoamericana},
    VOLUME = {38},
      YEAR = {2022},
    NUMBER = {4},
     PAGES = {1349--1370},
      ISSN = {0213-2230,2235-0616},
   MRCLASS = {35Q53 (53C42)},
  MRNUMBER = {4445917},
       DOI = {10.4171/rmi/1352},
       URL = {https://doi.org/10.4171/rmi/1352},
}

@article{catenary-cylinders,
    AUTHOR = {Mart\'{\i}nez-Trivi\~{n}o, A. L. and dos Santos, J. P.},
     TITLE = {Uniqueness of the {$[\varphi, \overrightarrow e_3]$}-catenary
              cylinders by their asymptotic behaviour},
   JOURNAL = {J. Math. Anal. Appl.},
  FJOURNAL = {Journal of Mathematical Analysis and Applications},
    VOLUME = {514},
      YEAR = {2022},
    NUMBER = {2},
     PAGES = {Paper No. 126347, 17},
      ISSN = {0022-247X,1096-0813},
   MRCLASS = {53A10 (35B33)},
  MRNUMBER = {4429585},
MRREVIEWER = {Ergin\ Bayram},
       DOI = {10.1016/j.jmaa.2022.126347},
       URL = {https://doi.org/10.1016/j.jmaa.2022.126347},
}

@misc{half-space,
  author = {Mart\'{\i}nez-Trivi\~{n}o, A. L. and dos Santos, J. P. and Tinaglia, T.},
  title  = {Half-space type theorems for a class of weighted minimal surfaces in $\mathbb{R}^3$},
  year   = {2025},
  note   = {arXiv:2601.22092},
  url    = {https://arxiv.org/abs/2601.22092}
}

@article{weierstrass-martinez-trivino,
    AUTHOR = {Mart\'{\i}nez, Antonio and Mart\'{\i}nez-Trivi\~{n}o, A. L.},
     TITLE = {A {W}eierstrass type representation for translating solitons
              and singular minimal surfaces},
   JOURNAL = {J. Math. Anal. Appl.},
  FJOURNAL = {Journal of Mathematical Analysis and Applications},
    VOLUME = {516},
      YEAR = {2022},
    NUMBER = {2},
     PAGES = {Paper No. 126528, 11},
      ISSN = {0022-247X,1096-0813},
   MRCLASS = {53A10 (35C08)},
  MRNUMBER = {4457891},
MRREVIEWER = {Pradip\ Kumar},
       DOI = {10.1016/j.jmaa.2022.126528},
       URL = {https://doi.org/10.1016/j.jmaa.2022.126528},
}

@article{lopez2018invariant,
  title={Invariant singular minimal surfaces},
  author={L{\'o}pez, Rafael},
  journal={Annals of Global Analysis and Geometry},
  volume={53},
  number={4},
  pages={521--541},
  year={2018},
  publisher={Springer}
}

@book{khalil2002nonlinear,
  title={Nonlinear systems},
  author={Khalil, Hassan K and Grizzle, Jessy W},
  volume={3},
  year={2002},
  publisher={Prentice hall Upper Saddle River, NJ}
}

@article{dierkes2021symmetric,
  title={Symmetric solutions of the singular minimal surface equation},
  author={Dierkes, Ulrich and Groh, Nico},
  journal={Annals of Global Analysis and Geometry},
  volume={60},
  number={2},
  pages={431--453},
  year={2021},
  publisher={Springer}
}

@article{lopez2019dirichlet,
  title={The Dirichlet problem for the $\alpha$-singular minimal surface equation},
  author={L{\'o}pez, Rafael},
  journal={Archiv der Mathematik},
  volume={112},
  number={2},
  pages={213--222},
  year={2019},
  publisher={Springer}
}

@article{martin2015topology,
  title={On the topology of translating solitons of the mean curvature flow},
  author={Mart{\'\i}n, Francisco and Savas-Halilaj, Andreas and Smoczyk, Knut},
  journal={Calculus of Variations and Partial Differential Equations},
  volume={54},
  number={3},
  pages={2853--2882},
  year={2015},
  publisher={Springer}
}

@misc{tinaglia-pipoli,
  author = {Pipoli, G. and dos Santos, J. P. and Tinaglia, T.},
  title  = {On the geometry of the asymptotic boundary of translators in $\mathbb{H}^2 \times \mathbb{R}$},
  year   = {2025},
  note   = {arXiv:2505.21083},
  url    = {https://arxiv.org/abs/2505.21083}
}

\end{document}